\documentclass[graybox]{svmult}
\pdfoutput=1

\usepackage{amsmath,amssymb,etex,ifxetex}
\usepackage{type1cm}
\usepackage{graphicx}
\usepackage{multicol}
\usepackage[bottom]{footmisc}
\usepackage{newtx}

\usepackage[utf8]{inputenc}
\usepackage{mathtools,ragged2e}
\usepackage{tikz}
\usetikzlibrary{calc,shapes,shapes.callouts,shapes.arrows,patterns,fit,backgrounds,decorations.pathmorphing}
\usepackage{booktabs,framed}
\usepackage[all]{xy}
\usepackage[protrusion=true,expansion=true]{microtype}
\usepackage{xspace}
\usepackage{marginnote,animate,ocgx2}

\usepackage{natbib}

\usepackage[{colorlinks,citecolor=[rgb]{0,0.1,0}}]{hyperref}

\graphicspath{{images/}}

\usepackage{pifont}
\newcommand{\cmark}{\ding{51}}
\newcommand{\xmark}{\ding{55}}
\definecolor{mypurple}{RGB}{150,0,255}

\usepackage{hyperref}

\newenvironment{indentblock}{%
  \list{}{\leftmargin\leftmargin}%
  \item\relax
}{%
  \endlist
}

\renewcommand{\E}{\mathcal{E}}

\newcommand{\BB}{\mathbb{B}}
\newcommand{\NN}{\mathbb{N}}
\newcommand{\QQ}{\mathbb{Q}}
\newcommand{\TT}{\mathbb{T}}
\newcommand{\RR}{\mathbb{R}}
\newcommand{\CC}{\mathbb{C}}

\newcommand{\defeq}{\vcentcolon=}

\DeclareMathOperator{\Sh}{Sh}
\DeclareMathOperator{\Zar}{Zar}

\newcommand{\Set}{\mathrm{Set}}
\newcommand{\Eff}{\mathrm{Ef{}f}}
\renewcommand{\_}{\mathpunct{.}\,}
\newcommand{\effective}{ef{}fective\xspace}
\newcommand{\effectively}{ef{}fectively\xspace}
\newcommand{\?}{\,{:}\,}
\newcommand{\realizes}{\Vdash}
\newcommand{\notnot}{\emph{not~not}\xspace}
\newcommand{\seq}[1]{\mathrel{\vdash\!\!\!_{#1}}}

\renewcommand{\paragraph}[1]{\noindent\textbf{#1.}}

\begin{document}

\addtocounter{chapter}{3}

\title{Exploring mathematical objects from custom-tailored mathematical universes}
\author{Ingo Blechschmidt}
\institute{Ingo Blechschmidt \at Universität Augsburg, Institut für Mathematik,
Universitätsstr. 14,
86159 Augsburg, Germany, \email{ingo.blechschmidt@math.uni-augsburg.de}}

\maketitle

\abstract{
  Toposes can be pictured as mathematical universes. Besides the standard topos,
  in which most of mathematics unfolds, there is a colorful host of alternate
  toposes in which mathematics plays out slightly
  differently.
  For instance, there are toposes in which the axiom of choice and the
  intermediate value theorem from undergraduate calculus fail.
  The purpose of this contribution is to give a glimpse of the toposophic
  landscape, presenting several specific toposes and exploring their peculiar
  properties, and to explicate how toposes provide distinct lenses through
  which the usual mathematical objects of the standard topos can be viewed.
}
\keywords{topos theory, realism debate, well-adapted language,
  constructive mathematics}
\medskip

\begin{figure}[b]
  \centering
  \tikzstyle{topos} = [draw=mypurple, very thick, rectangle, rounded corners, inner sep=5pt, inner ysep=10pt]
  \tikzstyle{title} = [fill=mypurple, text=white]


\newcommand{\setisprime}[1]{
  \ifnum#1=1 \gdef\isprime{0} \else \gdef\isprime{1} \fi
  \foreach \sip in {2, 3,5,...,#1} {
    \pgfmathparse{\sip*\sip>#1? 1:0}
    \ifthenelse{\pgfmathresult=1}{
      \breakforeach
    }{
      \pgfmathparse{Mod(#1,\sip)==0? 1:0}
      \ifthenelse{\pgfmathresult=1}{
        \gdef\isprime{0}
        \breakforeach
      }{}
    }
  }
}

\newcommand{\setxy}[1]{
  \pgfmathtruncatemacro{\x}{Mod(#1-1,\cols)}
  \pgfmathtruncatemacro{\y}{(#1-1) / \cols}
  \pgfmathtruncatemacro{\y}{\cols - 1 - \y}
  \pgfmathparse{2.5*(\x+.5)}\let\x\pgfmathresult
  \pgfmathparse{2.5*(\y+.5)}\let\y\pgfmathresult
}

\newcommand{\numlabel}[2]{
  \setxy{\n}
  \node[fill=none, text=black] at (\x,\y) {#2};
}

\newcommand{\drawpolygon}[2]{
  \setxy{#1}
  \ifthenelse{#2>1}{ 
    \ifthenelse{#2<30}{ 
      \filldraw (\x,\y) +(90:1)
      \foreach \drawi in {1,...,#2} {-- +(\drawi/#2*360+90:1)} -- cycle;
    }{ 
      \filldraw (\x,\y) circle(1);
    }
  }{}
}

\newcommand{\setpolygoncolor}[1]{
  \gdef\polycolor{black}
  \ifnum#1=2\gdef\polycolor{black!50!white}\fi
  \ifnum#1=3\gdef\polycolor{yellow!95!red}\fi
  \ifnum#1=5\gdef\polycolor{yellow!0!red}\fi
  \ifnum#1=7\gdef\polycolor{blue!75!green}\fi
  \ifnum#1=11\gdef\polycolor{blue!70!red}\fi
  \ifnum#1=13\gdef\polycolor{blue!40!red}\fi
  \ifnum#1=17\gdef\polycolor{green!50!blue}\fi
  \ifnum#1=19\gdef\polycolor{green!80!black}\fi
  \ifnum#1=23\gdef\polycolor{green!50!red}\fi
  \ifnum#1=29\gdef\polycolor{yellow!50!black}\fi
  \ifnum#1=31\gdef\polycolor{orange!50!black}\fi
  \ifnum#1=37\gdef\polycolor{red!50!black}\fi
  \ifnum#1=41\gdef\polycolor{purple!50!black}\fi
  \ifnum#1=43\gdef\polycolor{blue!50!black}\fi
  \ifnum#1=47\gdef\polycolor{green!50!black}\fi
  \ifnum#1=53\gdef\polycolor{white!50!black}\fi
  \ifnum#1=59\gdef\polycolor{white!50!black}\fi
  \ifnum#1=61\gdef\polycolor{white!50!black}\fi
  \ifnum#1=67\gdef\polycolor{white!50!black}\fi
}

\newcommand{\sieve}[2]{
  \def\cols{#1}
  \def\rows{#2}
  \begin{tikzpicture}[scale=.5]
  \pgfmathtruncatemacro{\nmax}{\rows * \cols}

  \foreach \n in {1,...,\nmax} {
    \begin{scope}[fill=gray, fill opacity=.05,
                  draw=gray, draw opacity=.10,
                  line width=4]
      \drawpolygon{\n}{\n}
    \end{scope}
    \setisprime{\n}
    \ifthenelse{\isprime=1}{
      \numlabel{\n}{\bf\n}
    }{
      \def\startintensity{.33}
      \def\incrintensity{.10}
      \def\intensity{\startintensity}

      \def\m{\n}
      \pgfmathtruncatemacro{\i}{\m / 2}

      \whiledo{\m>1}{
        \setisprime{\i}
        \pgfmathparse{Mod(\m,\i)==0? 1:0}
        \ifthenelse{\pgfmathresult=1\and\isprime=1}{
          \setpolygoncolor{\i}
          \begin{scope}[fill=\polycolor, fill opacity=\intensity,
                        draw=\polycolor!85!black, draw opacity=\intensity,
                        line width=\intensity*1.5]
            \drawpolygon{\n}{\i}
          \end{scope}
          \pgfmathtruncatemacro{\m}{\m / \i}
          \pgfmathparse{\intensity + \incrintensity}\let\intensity\pgfmathresult
        }{
          \pgfmathtruncatemacro{\i}{\i - 1}
          \def\intensity{\startintensity}
        }
      }
      \begin{scope}[text=black, text opacity=.5]
        \numlabel{\n}{\scriptsize\n}
      \end{scope}
    }
  }

  \end{tikzpicture}
}

\newcommand{\fakesieve}[2]{
  \def\cols{#1}
  \def\rows{#2}
  \begin{tikzpicture}[scale=.5,opacity=0]
  \pgfmathtruncatemacro{\nmax}{\rows * \cols}

  \foreach \n in {1,...,\nmax} {
    \begin{scope}[fill=gray,
                  draw=gray,
                  line width=4]
      \drawpolygon{\n}{\n}
    \end{scope}
    \setisprime{\n}
    \ifthenelse{\isprime=1}{
      \numlabel{\n}{\bf\n}
    }{
      \def\startintensity{.33}
      \def\incrintensity{.10}
      \def\intensity{\startintensity}

      \def\m{\n}
      \pgfmathtruncatemacro{\i}{\m / 2}

      \whiledo{\m>1}{
        \setisprime{\i}
        \pgfmathparse{Mod(\m,\i)==0? 1:0}
        \ifthenelse{\pgfmathresult=1\and\isprime=1}{
          \setpolygoncolor{\i}
          \begin{scope}[fill=\polycolor,
                        draw=\polycolor!85!black,
                        line width=\intensity*1.5]
            \drawpolygon{\n}{\i}
          \end{scope}
          \pgfmathtruncatemacro{\m}{\m / \i}
          \pgfmathparse{\intensity + \incrintensity}\let\intensity\pgfmathresult
        }{
          \pgfmathtruncatemacro{\i}{\i - 1}
          \def\intensity{\startintensity}
        }
      }
      \begin{scope}[text=black]
        \numlabel{\n}{\scriptsize\n}
      \end{scope}
    }
  }

  \end{tikzpicture}
}

  \newcommand{\drawbox}[4]{
    \node[topos, #4] [fit = #3] (#1) {};
    \node[title] at (#1.north) {#2};
  }

  \newcommand{\muchstuff}{
    \includegraphics[height=3em]{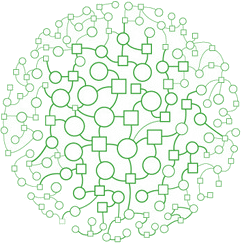}
    \scalebox{0.5}{\sieve{14}{2}}
  }

  \newcommand{\muchstuffplaceholder}{
    \includegraphics[height=3em]{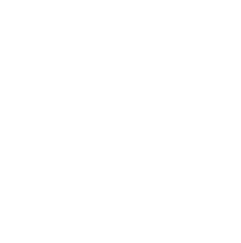}
    \scalebox{0.5}{\fakesieve{14}{2}}
  }

  \newcommand{\fewstuff}{
    \includegraphics[height=3em]{filmat}
    \scalebox{0.5}{\sieve{7}{2}}
  }

  \newcommand{\threeblobs}{
    \colorbox{mypurple}{\ \ }\quad
    \colorbox{mypurple}{\ \ }\quad
    \colorbox{mypurple}{\ \ }
  }

  \begin{tikzpicture}
    \node (objs-set0) at (0,0) {
      \muchstuff
    };
    \node[scale=0.4] (objs-set1) at (-3.5,-2.5) {
      \fewstuff
    };
    \node[scale=0.4] (objs-eff1) at (3.5,-2.5) {
      \fewstuff
    };
    \node[scale=0.4] (objs-sh1)  at (0,-2.5) {
      \fewstuff
    };

    \node (prop-set1) [below of=objs-set1, align=left] {
      The usual laws \\
      of logic hold.
    };

    \node (prop-eff1) [below of=objs-eff1, align=left] {
      Every function \\
      is computable.
    };

    \node (prop-sh1) [below of=objs-sh1, align=left] {
      The axiom of \\
      choice fails.
    };

    \node (more-eff1) [below of=prop-eff1] {
      \threeblobs
    };
    \node (more-sh1)  [below of=prop-sh1] {
      \threeblobs
    };
    \node (more-set1) [below of=prop-set1] {
      \threeblobs
    };

    \drawbox{set1}{$\mathrm{Set}$}{(objs-set1) (prop-set1) (more-set1)}{}
    \drawbox{eff1}{Ef{}f}{(objs-eff1) (prop-eff1) (more-eff1)}{tape}
    \drawbox{sh1}{$\mathrm{Sh}\, X$}{(objs-sh1) (prop-sh1) (more-sh1)}{draw=none}
    \def\R{8pt}
    \begin{pgfonlayer}{background}
      \draw[decoration={bumps,segment length=8pt}, decorate, very thick, draw=mypurple]
        ($(sh1.south west) + (\R, 0)$) arc(270:180:\R) --
        ($(sh1.north west) + (0, -\R)$) arc(180:90:\R) --
        ($(sh1.north east) + (-\R, 0)$) arc(90:0:\R) --
        ($(sh1.south east) + (0, \R)$) arc(0:-90:\R) --
        cycle;
    \end{pgfonlayer}
    \drawbox{set0}{$\mathrm{Set}$}{(objs-set0) (set1) (eff1) (sh1)}{}
  \end{tikzpicture}

  \caption{\label{fig:landscape}A glimpse of the toposophic landscape,
  displaying alongside the standard topos~$\Set$ two further toposes.}
\end{figure}
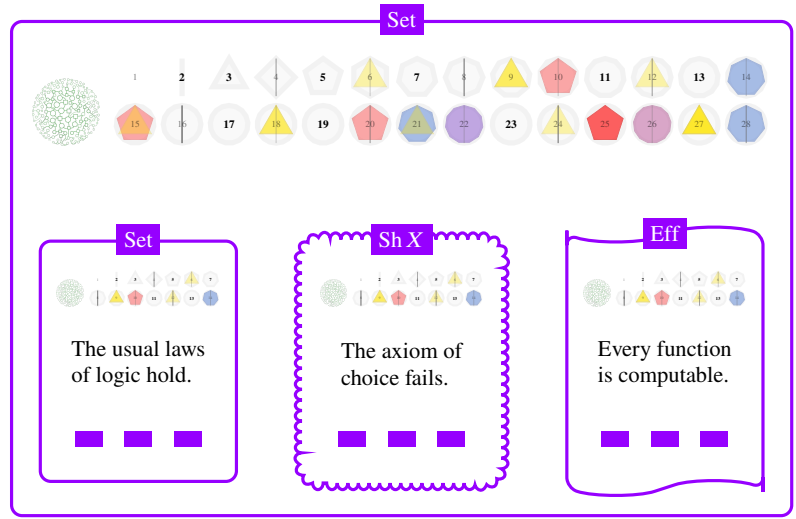

\noindent
Toposes can be pictured as mathematical universes in which we can do
mathematics. Most mathematicians spend all their professional life in just a
single topos, the so-called \emph{standard topos}. However, besides the
standard topos, there is a colorful host of alternate toposes which are just as
worthy of mathematical study and in which mathematics plays out slightly
differently (Figure~\ref{fig:landscape}).

For instance, there are toposes in which the axiom of choice and the
intermediate value theorem from undergraduate calculus fail, toposes in which
any function~$\RR \to \RR$ is continuous and toposes in which infinitesimal
numbers exist.

The purpose of this contribution is twofold.
\begin{enumerate}
\item We give a glimpse of the toposophic landscape, presenting several
specific toposes and exploring their peculiar properties.

\item We explicate how toposes provide distinct lenses through which the
usual mathematical objects of the standard topos can be viewed.
\end{enumerate}

Viewed through such a lens, a given mathematical object can have different
properties than when viewed normally. In particular, it can have
better properties for the purposes of specific applications, especially if
the topos is custom-tailored to the object in question. This change of
perspective has been used in mathematical practice.
To give just a taste of what
is possible, through the lens provided by an appropriate topos, any given ring
can look like a field and hence mathematical techniques for fields also apply,
through the lens, to rings.

We argue that toposes and specifically the change in perspective provided by
toposes are ripe for philosophical analysis. In particular, there are the
following connections with topics in the philosophy of mathematics:
\begin{enumerate}
\item Toposes enrich the realism/anti-realism debate in that they paint the larger
picture that the platonic heaven of mathematical objects is not unique: besides
the standard heaven of the standard topos, we can fathom the alternate
heavens of all other toposes, all embedded in a second-order heaven.
\item To some extent, the mathematical landscape depends on the commonly agreed-upon rules of
mathematics. These are not entirely absolute; for instance, it is conceivable
that from the foundational crisis Brouwer's intuitionism would have emerged as
the main school of thought and that we would now all reject the law of excluded
middle. Toposes allow us to explore alternatives to how history
has played out.
\item Mathematics is not only about studying mathematical objects, but also
about studying the relations between mathematical objects. The distinct view
on mathematical objects provided by any topos uncovers relations which
otherwise remain hidden.\footnote{The research
program put forward by~\cite{caramello:tst} provides a further topos-theoretic way for
uncovering hidden relations, though not between objects but between mathematical
theories. \textit{\textcolor{red}{Note to editor: please insert cross reference to Olivia's chapter of
this book.}}}
\item In some cases, a mathematical relation can be expressed quite succinctly
using the language of a specific topos and not so succinctly using the language
of the standard topos. This phenomenon showcases the importance of
\emph{appropriate language}.
\item Toposes provide new impetus to study constructive mathematics and
intuitionistic logic, in particular also to restrict to intuitionistic
logic on the meta level and to consider the idea that the platonic heaven might
be governed by intuitionistic logic.
\end{enumerate}
We invite further research on these connections.

We intend this contribution to be self-contained and do not assume familiarity with
topos theory or category theory, having a diverse readership of people interested in philosophy of
mathematics in mind. However, to make this text more substantial to
categorically-inclined readers, some categorical definitions are included.
These definitions can be skipped without impacting the main message of this contribution.

Readers who would like to learn more details
are directed to the survey of category theory by~\cite{sep:category-theory}
and to a gentle introduction to topos theory
by~\cite{leinster:introduction}. Standard references for the internal
language of toposes include~\cite[Chapter~VI]{moerdijk-maclane:sheaves-logic},
\cite[Chapter~14]{goldblatt:topoi},
\cite{caramello:preliminaries}, \cite{streicher:ctcl}, \cite{shulman:categorical-logic},
\cite[Chapter~6]{borceux:handbook3} and~\cite[Part~D]{johnstone:elephant}.

\bigskip
\paragraph{Other aspects of toposes} This note focuses on just a single
aspect of toposes, the view of toposes as alternate mathematical universes.
This aspect is not the only one, nor did it historically come first.

Toposes were originally conceived by Grothendieck in the early 1960s for the
needs of algebraic geometry, as a general framework for constructing and
studying invariants in classical and new geometric contexts, and it is in that
subject that toposes saw their deepest applications. The proof of Fermat's Last
Theorem is probably the most prominent such application, crucially resting on
the cohomology and homotopy invariants provided by toposes.

In the seminal work introducing
toposes by~\cite{artin-grothendieck-verdier:topos}, toposes are viewed as
generalized kinds of spaces. Every topological space~$X$ gives rise to a topos,
the \emph{topos of sheaves over~$X$}, and every continuous map gives rise to a
\emph{geometric morphism} between the induced sheaf toposes, but not every
topos is of this form. While the open sets of a topological space are required
to be parts of the space, the opens of toposes are not; and while for open
subsets~$U$ and~$V$ there is only a truth value as to whether~$U$ is contained
in~$V$, in a general topos there can be many distinct ways how an open is
contained in another one. This additional flexibility is required in situations
where honest open subsets are rare, such as when studying the étale cohomology
of a scheme as in~\cite{milne:etale}.

That toposes could also be regarded as mathematical universes was realized only
later, by Bill Lawvere and Myles Tierney at the end of the 1960s. They abstracted some of
the most important categorical properties of Grothendieck's toposes into what
is now known as the definition of an \emph{elementary topos}. Elementary
toposes are considerably more general and less tied to geometry than the
original toposes. The theory of elementary toposes has a substantially
different, logical flavor, not least because a different notion of morphism
plays an important role. To help disambiguate, there is a trend to rename
elementary toposes to \emph{logoses}, but this text still follows the standard
convention.

A further perspective on toposes emerged in the early 1970s with the discovery
that toposes can be regarded as embodiments of a certain kind of first-order
theories, the \emph{geometric theories} briefly discussed on
page~\pageref{item:classifying-topos}. The so-called \emph{classifying toposes}
link geometrical and logical aspects and are fundamental to Olivia Caramello's
bridge-building program set out in~\cite{caramello:tst}. Geometrically, the classifying
topos of a geometric theory~$\TT$ can be regarded as the generalized space of
models of~$\TT$; this idea is due to~\cite{hakim:relative-schemes}, though she did not cast her discovery in
this language. Logically, the classifying topos of~$\TT$ can be regarded as a
particular mathematical universe containing the \emph{generic~$\TT$-model}, a
model which has exactly those properties which are shared by all models.

Yet more views on toposes are fruitfully employed --
\cite[pages~vii--viii]{johnstone:elephant} lists ten more -- but we shall not review
them here. A historical survey was compiled by~\cite{mclarty:history}.

\bigskip
\paragraph{Acknowledgments} We are grateful to Andrej Bauer, Martin
Brandenburg, Sina Hazratpour, Matthias Hutzler, Marc Nieper-Wißkirchen and
Alexander Oldenziel for invaluable discussions shaping this work, to Moritz
Laudahn, Matthias Hutzler and Milan Zerbin for their careful criticism of
earlier drafts and to Todd Lehman for the code for parts of
Figure~\ref{fig:landscape}. This note also profited substantially from the
thorough review by three anonymous referees, whose efforts are much
appreciated and gladly acknowledged. We thank the organizers and all participants of the
Mussomeli conference of the Italian Network for the Philosophy in Mathematics,
where this work was presented, for creating an exceptionally beautiful meeting. In
particular, we are grateful to Neil Barton, Danielle Macbeth, Gianluigi Oliveri
and Lorenzo Rossi for valuable comments.


\section{Toposes as alternate mathematical universes}

A topos is a certain kind of \emph{category}, containing objects and
morphisms between those objects. The precise definition is recorded here only for
reference. Appreciating it requires some amount of category theory, but, as will be demonstrated in
the following sections, exploring the mathematical universe of a given topos
does not.

\begin{definition}A \emph{topos} is a category which has all finite
limits, is cartesian closed, has a subobject classifier and contains a natural
numbers object.\footnote{More precisely, this is the definition of an
\emph{elementary topos with a natural numbers object}. Since this definition is
less tied to geometry than Grothendieck's (as categories of sheaves over
sites), there is a trend to call these toposes \emph{logoses}. However, that
term also has other uses.}\end{definition}

Put briefly, these axioms state that a topos should share several categorical
properties with the category of sets; they ensure that each topos contains its
own versions of familiar mathematical objects such as natural numbers, real
numbers, groups and manifolds, and is closed under the usual constructions
such as cartesian products or quotients.
The prototypical topos is the standard topos:

\begin{definition}The \emph{standard topos}~$\Set$ is the category which has all
sets as its objects and all maps between sets as morphisms.\end{definition}

Given a topos~$\E$, we write~``$\E \models \varphi$'' to denote that a
mathematical statement~$\varphi$ \emph{holds in~$\E$}. The meaning of~``$\E \models
\varphi$'' is defined by recursion on the structure of~$\varphi$ following the
so-called \emph{Kripke--Joyal translation rules}. For instance, the rules for
translating conjunction and falsity read
\[ \begin{array}{lcl}
  \E \models (\alpha \wedge \beta) &\qquad\text{iff}\qquad&
  \E \models \alpha \quad\text{and}\quad \E \models \beta, \\[0.2em]
  \E \models \bot &\qquad\text{iff}\qquad&
  \text{$\E$ is the trivial topos}.
\end{array} \]
The remaining translation rules are more
involved, as detailed by~\cite[Section~VI.7]{moerdijk-maclane:sheaves-logic}; we do not list them here for
the case of a general topos~$\E$, but we will state them in the next sections
for several specific toposes. We refer to~``$\E \models \varphi$'' also as the
``external meaning of the internal statement~$\varphi$''.

In the definition of~$\E \models \varphi$, the statement~$\varphi$ can be any
statement in the language of a general version of higher-order predicate
calculus with dependent types, with a base type for each object of~$\E$ and
with a constant of type~$X$ for each morphism~$1 \to X$ in~$\E$. In practice almost any mathematical statement
can be interpreted in a given topos.\footnote{The main exceptions are
statements from set theory, which typically make substantial use of a global
membership predicate~``$\in$''. Toposes only support a typed \emph{local}
membership predicate, where we may write~``$x \in A$'' only in the context of
some fixed type~$M$ such that~$x$ is of type~$M$ and~$A$ is of
type~$P(M)$, the power type of~$M$. We refer
to~\cite{fourman:sheaf-models,streicher:forcizf,awodey-butz-simpson-streicher:bist}
for ways around this restriction.} We refrain from giving a precise definition
of the language here, but refer to the references~\cite[Section~7]{shulman:stack}
and~\cite[Section~VI.7]{moerdijk-maclane:sheaves-logic} for details.

It is by the Kripke--Joyal translation rules that we can access the alternate
universe of a topos. In the special case of the standard topos~$\Set$, the
definition of~``$\Set \models \varphi$'' unfolds to~$\varphi$ for any
statement~$\varphi$. Hence a statement holds in the standard topos if and only
if it holds in the usual mathematical sense.


\subsection{The logic of toposes}\label{sect:logic-of-toposes}
By their definition as special kinds of
categories, toposes are merely algebraic structures not unlike groups or vector
spaces. Hence we need to argue why we picture toposes as mathematical universes
while we do not elevate other kinds of algebraic structures in the same way.
For us, this usage is justified by the following metatheorem:

\begin{theorem}\label{thm:reasoning}Let~$\E$ be a topos and let~$\varphi$ be a
statement such that~$\E \models \varphi$. If~$\varphi$ intuitionistically
entails a further statement~$\psi$ (that is, if it is provable in
intuitionistic logic that~$\varphi$ entails~$\psi$), then~$\E \models
\psi$.\end{theorem}

This metatheorem allows us to \emph{reason} in toposes. When first exploring a
new topos~$\E$, we need to employ the Kripke--Joyal translation rules each time
we want to check whether a statement holds in~$\E$. But as soon as we
have amassed a stock of statements known to be true in~$\E$, we can find more
by deducing their logical consequences.

For instance, in any topos where the statement ``any map~$\RR \to \RR$ is
continuous'' is true, also the statement ``any map~$\RR \to \RR^2$ is
continuous'' is, since there is an intuitionistic proof that a map into a
higher-dimensional Euclidean space is continuous if its individual components
are.

The only caveat of Theorem~\ref{thm:reasoning} is that toposes generally only
support intuitionistic reasoning and not the full power of the ordinary
\emph{classical reasoning}. That is, within most toposes, the law of excluded
middle ($\varphi \vee \neg\varphi$) and the law of double negation elimination
($\neg\neg\varphi \Rightarrow \varphi$) are not available. It is intuitionistic
logic and not classical logic which is the common denominator of all toposes;
we cannot generally argue by contradiction in a topos.

While it may appear that these two laws pervade any mathematical theory, in
fact a substantial amount of mathematics can be developed intuitionistically
(see for
instance~\cite{mines-richman-ruitenburg:constructive-algebra,lombardi-quitte:constructive-algebra}
for constructive algebra,~\cite{bishop-bridges:constructive-analysis} for
constructive analysis
and~\cite{bauer:int-mathematics,bauer:video,melikhov:intuitionistic-logic} for
accessible surveys on appreciating intuitionistic logic) and hence the alternate universes provided by toposes
cannot be too strange: In any topos, there are infinitely many prime numbers,
the square root of two is not rational, the fundamental theorem of Galois
theory holds and the powerset of the naturals is uncountable.

That said, intuitionistic logic still allows for a considerable amount of
freedom, and in many toposes statements are true which are baffling if one has
only received training in mathematics based on classical logic. For instance, on first sight it
looks like the sign function
\[ \operatorname{sgn} : \RR \longrightarrow \RR,\ x \longmapsto \begin{cases}
  -1, & \text{if $x < 0$,} \\
  0, & \text{if $x = 0$,} \\
  1, & \text{if $x > 0$,}
\end{cases} \]
is an obvious counterexample to the statement ``any map~$\RR \to \RR$ is
continuous''. However, a closer inspection reveals that the sign function
cannot be proven to be a total function~$\RR \to \RR$ if only intuitionistic
logic is available. The domain of the sign
function is the subset~$\{ x \in \RR \,|\, x < 0 \vee x = 0 \vee x > 0 \}
\subseteq \RR$, and in intuitionistic logic this subset cannot be shown to
coincide with~$\RR$.

Sections~\ref{sect:effective-topos} to~\ref{sect:smooth} present several
examples for such anti-classical statements and explain how to make sense of
them. There are also toposes which are closer to the standard topos and do not validate such
anti-classical statements:

\begin{definition}A topos~$\E$ is \emph{boolean} if and only if the laws of classical
logic are true in~$\E$.\end{definition}

Since exactly those statements hold in the standard topos which hold on the
meta level, the standard topos is boolean if and only if, as is commonly supposed, the laws
of classical logic hold on the meta level. Most toposes of interest are not
boolean, irrespective of one's philosophical commitments about the meta level,
and conversely some toposes are boolean even if classical logic is not
available on the meta level.

\begin{remark}The axiom of choice (which is strictly speaking not part of
classical logic, but of classical set theory) is also not available in most
toposes. By \emph{Diaconescu's theorem}, the axiom of choice implies the law of
excluded middle in presence of other axioms which are available in any topos.
\end{remark}



At this point in the text, all prerequisites for exploring toposes have been
introduced. The reader who wishes to develop, by explicit examples, intuition
for working internally to toposes is invited to skip ahead to
Section~\ref{sect:effective-topos}.

\subsection{Relation to models of set theory} In set theory, philosophy and
logic, models of set theories are studied. These are structures~$(M,\in)$
validating the axioms of some set theory such as Zermelo--Fraenkel set theory
with choice~\textsc{zfc}, and they can be pictured as ``universes in which we
can do mathematics'' in much the same way as toposes.

In fact, to any model~$(M,\in)$ of a set theory such as~\textsc{zf}
or~\textsc{zfc}, there is a topos~$\Set_M$ such that a statement holds
in~$\Set_M$ if and only if it holds in~$M$.\footnote{The topos~$\Set_M$ can be
described as follows: Its objects are the elements of~$M$, that is the entities
which~$M$ believes to be sets, and its morphisms are those entities which~$M$
believes to be maps. The topos~$\Set_M$ validates the axioms of the structural
set theory \textsc{etcs}, see~\cite{mclarty:structuralism,marquis:foundations,barton-friedman:structures}, and models are isomorphic if and only if their
associated toposes are equivalent as
categories, see~\cite[Section~VI.10]{moerdijk-maclane:sheaves-logic}.}

\begin{example}The topos~$\Set_V$ associated to the universe~$V$ of all sets (if
this structure is available in one's chosen ontology) coincides with the
standard topos~$\Set$.\end{example}

In set theory, we use forcing and other techniques to construct new
models of set theory from given ones, thereby exploring the set-theoretic
multiverse. There are similar techniques available for constructing new toposes
from given ones, and some of these correspond to the techniques from set
theory.

However, there are also important differences between the notion of mathematical
universes as provided by toposes and as provided by models of set theory, both
regarding the subject matter and the reasons for why we are interested in them.

Firstly, toposes are more general than models of set theory. Every model of set
theory gives rise to a topos, but not every topos is induced in this way from a model of set
theory. Unlike models of~\textsc{zfc}, most toposes do not validate the law of excluded middle, much
less so the axiom of choice.

Secondly, there is a shift in emphasis. An important philosophical objective
for studying models of set theory is to explore which notions of sets are
coherent: Does the cardinality of the reals need to be the cardinal directly
succeeding~$\aleph_0$, the cardinality of the naturals? No, there are models of
set theory in which the continuum hypothesis fails. Do non-measurable sets of
reals need to exist? No, in models of~$\textsc{zf}+\textsc{ad}$,
Zermelo--Fraenkel set theory plus the axiom of determinacy, it is a theorem that
every subset of~$\RR^n$ is Lebesgue-measurable. Can the axiom of choice be
added to the axioms of~\textsc{zf} without causing inconsistency? Yes, if~$M$
is a model of~\textsc{zf} then~$L^M$, the structure of the constructible sets of~$M$, forms a
model of~\textsc{zfc}.

Toposes can be used for similar such purposes, and indeed have been,
especially to explore the various intuitionistic notions of sets. However, an important
aspect of topos theory is that toposes are used to explore the \emph{standard}
mathematical universe: truth in the \effective topos tells us what is
computable; truth in sheaf toposes tells us what is true locally; toposes
adapted to synthetic differential geometry can be used to rigorously work with
infinitesimals. All of these examples will be presented in more detail in the
next sections.

In a sense which can be made precise, toposes allow us to study the usual
objects of mathematics from a different point of view -- one such view for
every topos -- and it is a beautiful and intriguing fact that with the sole
exception of the law of excluded middle, the laws of logic apply to
mathematical objects also when viewed through the lens of a specific topos.

\subsection{A glimpse of the toposophic landscape}
There is a proper class of toposes. Figure~\ref{fig:landscape} depicts three
toposes side by side: the standard topos, a sheaf topos and the \effective
topos. Each of these toposes tells a different story of mathematics, and any
topos which is not the standard topos invites us to ponder alternative ways how
mathematics could unfold.

Some of the most prominent toposes are the following.

\begin{enumerate}
\item The \emph{trivial topos}. In the trivial topos, any statement whatsoever
is true. The trivial topos is not interesting on its own, but its existence
streamlines the theory and it can be an interesting question whether a given
topos coincides with the trivial topos.
\item $\Set$, the \emph{standard topos}. A statement is true in~$\Set$ iff it is true
in the ordinary mathematical sense.
\item $\Set_M$, the topos associated to any model~$(M,\in)$ of~\textsc{zf}.
\item $\Set^W$, the category of functors~$(W,\leq) \to \Set$ associated to
any Kripke model~$(W,\leq)$. A statement is true in this topos iff it is valid
with respect to the ordinary Kripke semantics of~$(W,\leq)$. This example shows
that the Kripke--Joyal semantics of toposes generalizes the more familiar
Kripke semantics.
\item $\Eff$, the \emph{\effective topos}. A statement is true in~$\Eff$ iff
if it has a \emph{computable witness} as detailed in
Section~\ref{sect:effective-topos}. In~$\Eff$, any function~$\NN \to \NN$ is
computable, any function~$\RR \to \RR$ is continuous and the countable axiom of
choice holds (even if it does not on the meta level).
\item $\Sh(X)$, the \emph{topos of sheaves} over any space~$X$. A statement is true
in~$\Sh(X)$ iff it holds \emph{locally on~$X$}, as detailed in
Section~\ref{sect:sheaf-toposes}. For most choices of~$X$, the
axiom of choice and the intermediate value theorem fail in~$\Sh(X)$, and this
failure is for geometric reasons.
\item $\Zar(A)$, the \emph{Zariski topos} of a ring~$A$ presented in
Section~\ref{sect:smooth}. This topos contains a mirror
image of~$A$ which is a field, even if~$A$ is not.
\item $\mathrm{Bohr}(A)$, the \emph{Bohr topos} associated to a noncommutative
C\textsuperscript{*}\kern-.1ex-algebra~$A$. This topos contains a mirror image
of~$A$ which is commutative. In this sense, quantum mechanical systems (which are
described by noncommutative C\textsuperscript{*}\kern-.1ex-algebras) can be
regarded as classical mechanical systems (which are described by commutative
algebras). Details are described by~\cite{butterfield-hamilton-isham:bohr,heunen-landsman-spitters:aqt}.
\item\label{item:classifying-topos}$\Set[\TT]$, the \emph{classifying topos} of a geometric theory~$\TT$.\footnote{A
geometric theory is a theory in many-sorted first-order logic whose axioms can be put as \emph{geometric
sequents}, sequents of the form~$\varphi \seq{\vec x} \psi$ where~$\varphi$
and~$\psi$ are geometric formulas (formulas built from equality and specified
relation symbols by the logical
connectives~${\top}\,{\bot}\,{\wedge}\,{\vee}\,{\exists}$ and by arbitrary
set-indexed disjunctions~$\bigvee$).} This topos contains the
\emph{generic~$\TT$-model}. For instance, the classifying topos of the theory
of groups contains the \emph{generic group}. Arguably it is this group
which we implicitly refer to when we utter the phrase ``Let~$G$ be a group.''.
The generic group has exactly those
properties which are shared by any group whatsoever.\footnote{More precisely,
this is only true for those properties which can be formulated as geometric
sequents. For arbitrary properties~$\varphi$, the statements~``the generic
group has property~$\varphi$'' and~``all groups have property~$\varphi$'' need
not be equivalent. This inbalance has mathematical applications and is explored
in~\cite{blechschmidt:nullstellensatz}.}
\item $T(\mathcal{L}_0)$, the \emph{free topos}. A statement is true in the free topos iff it is
intuitionistically provable. Lambek and Scott proposed that the free topos can
reconcile moderate platonism (because this topos has a certain universal
property which can be used to single it out among the plenitude of toposes),
moderate formalism (because it is constructed in a purely syntactic way) and
moderate logicism (because, as a topos, it supports an intuitionistic type
theory). Details are described by~\cite{lambek:incompatible,couture-lambek:reflections}.
\end{enumerate}

There are several constructions which produce new toposes from a given
topos~$\E$. A non-exhaustive list is the following.
\begin{enumerate}
\item Given an object~$X$ of~$\E$, the \emph{slice topos}~$\E/X$
contains the \emph{generic element}~$x_0$ of~$X$. This generic element can be
pictured as the element we implicitly refer to when we utter the
phrase~``Let~$x$ be an element of~$X$.''. A statement~$\varphi(x_0)$
about~$x_0$ is true in~$\E/X$ if and only if in~$\E$ the statement~$\forall x
\? X\_ \varphi(x)$ is true.

For instance, the topos~$\Set/\QQ$ contains the generic rational number~$x_0$.
Neither the statement~``$x_0$ is zero'' nor the statement~``$x_0$ is not zero''
hold in~$\Set/\QQ$, as it is neither the case that any rational number
in~$\Set$ is zero nor that any rational number in~$\Set$ is not zero. Like any
rational number, the number~$x_0$ can be written as a fraction~$\frac{a}{b}$.
Just as~$x_0$ itself, the numbers~$a$ and~$b$ are quite indetermined.
\item Given a statement~$\varphi$ (which may contain objects of~$\E$ as
parameters but which must be formalizable as a geometric sequent), there is a
largest subtopos of~$\E$ in which~$\varphi$ holds. This construction is useful
if neither~$\varphi$ nor~$\neg\varphi$ hold in~$\E$ and we want to
force~$\varphi$ to be true. If~$\E \models \neg\varphi$, then the resulting
topos is the trivial topos. (A subtopos is not simply a subcategory; rather, it
is more like a certain kind of quotient category. We do not give, and for the
purposes of this contribution do not need, further details.)
\item There is a ``smallest dense'' subtopos $\Sh_{\neg\neg}(\E)$. This topos
is always boolean, even if~$\E$ and the meta level are not. For a mathematician who employs
intuitionistic logic on their meta level, the nonconstructive results of their
classical colleagues do not appear to make sense in~$\Set$, but they hold
in~$\Sh_{\neg\neg}(\Set)$. If classical logic holds on the meta level,
then~$\Set$ and~$\Sh_{\neg\neg}(\Set)$ coincide.

The topos~$\Sh_{\neg\neg}(\E)$ is related to the \emph{double negation
translation} from classical logic into
intuitionistic logic: A statement holds in~$\Sh_{\neg\neg}(\E)$ if
and only if its translation holds
in~$\E$~\cite[Theorem~6.31]{blechschmidt:phd}.
\end{enumerate}

Toposes are still mathematical structures, and as long as we study
toposes within the usual setup of mathematics, our toposes are all part of the
standard topos. This is why Figure~\ref{fig:landscape} pictures the standard
topos twice, once as a particular topos next to others, and once as the universe
covering the entirety of our mathematical discourse.\footnote{There is a fine print to consider. Technically, if we
work within~\textsc{zf} or its intuitionistic cousin~\textsc{izf}, most toposes
of interest are proper classes, not sets. In particular~$\Set$ itself is a
proper class. Hence Figure~\ref{fig:landscape} should not be interpreted as
indicating that toposes are contained in~$\Set$ as objects, which most are not.
In this regard toposes are similar to class-sized inner models in set theory.
We believe that the vague statement ``our toposes are all part of the standard topos'' is still an
apt description of the situation. A possible formalization is (the trivial
observation that) ``our toposes are all indexed categories over~$\Set$''.}
The toposes which we can study in mathematics do not tell us
all possible stories how mathematics could unfold, only those which appear
coherent from the point of view of the standard topos, and the topos-theoretic
multiverse which we have access to is just a small part of an even larger
landscape.\footnote{This paragraph employs an overly narrow conception of
``mathematics'', focusing only on those mathematical worlds which form toposes
and for instance excluding any predicative flavors of mathematics
(Laura Crosilla's survey in~\cite{crosilla:predicativity} is an excellent introduction).
Toposes are impredicative in the sense that any object of a topos is required to have a
powerobject. A predicative cousin of toposes are the \emph{arithmetic
universes} introduced by Joyal which have recently been an important object of
consideration by Maietti and
Vickers, see~\cite{maietti:au,maietti-vickers:induction,vickers:sketches}.}

To obtain just a hint of how the true landscape looks like, we can study topos
theory from the inside of toposes; the resulting picture can look quite
different from the picture which emerges from within the standard topos.

For instance, from within the standard topos, we can write down the
construction which yields the standard topos and the construction which yields
the \effective topos~$\Eff$ and observe that the resulting toposes are not at
all equivalent: In~$\Eff$, any function~$\RR \to \RR$ is continuous
while~$\Set$ abounds with discontinuous functions (at least if we assume a
classical meta level). In contrast, if we carry out these two constructions
from within the \effective topos, we obtain toposes which are elementarily
equivalent. More precisely, for any statement~$\varphi$ of higher order
arithmetic,
\[ \Eff \models (\Set \models \varphi) \qquad\text{iff}\qquad\Eff \models
  (\Eff \models \varphi). \]
In this sense the construction which yields the \effective topos is
\emph{idempotent}~\cite[Section~3.8.3]{oosten:realizability}.

\begin{remark}The picture of a topos-theoretic multiverse is related to Hamkin's
multiverse view in set theory as put forward in~\cite{hamkins:multiverse}. In fact, the
topos-theoretic multiverse can be regarded as an extension of the set-theoretic
multiverse: While Hamkins proposes to embrace all models of set theory (not
necessarily all of them equally -- we might prefer some models over others), we propose to embrace all toposes (again
not necessarily all of them equally). As every model~$M$ of set theory
gives rise to a topos~$\Set_M$, the set-theoretic multiverse is contained in the
topos-theoretic one.

However, a central and intriguing feature of the multiverse view in set theory
has, as of yet, no counterpart in topos theory: namely a systematic study of
its modal logic with respect to various notions of relations between toposes.\end{remark}

\subsection{A syntactic account of toposes} We introduced toposes from a semantic
point of view. There is also a second, purely syntactic point of view on
toposes:
\begin{enumerate}
\item (semantic view) A topos is an alternate mathematical universe. Any topos
contains its own stock of mathematical objects. A ``transfer theorem'' relates
properties of those objects with properties of objects of the standard topos: A
statement~$\varphi$ about the objects of a topos~$\E$ holds in~$\E$ iff the
statement~``$\E \models \varphi$'' holds in the standard topos.
\item (syntactic view) A topos is merely an index to a syntactical translation
procedure. Any topos~$\E$ gives rise to a ``generalized modal operator'' which
turns a statement~$\varphi$ (about ordinary mathematical objects) into the
statement~``$\E \models \varphi$'' of the same kind (again about ordinary
mathematical objects).
\end{enumerate}

For instance, in the semantic view, the \effective topos is an alternative
universe which contains its own version of the natural numbers. These naturals
cannot be directly compared with the naturals of the standard topos, for they
live in distinct universes, but by the transfer theorem they are still linked in a
nontrivial way: For instance, the statement ``there are infinitely many primes
in~$\Eff$'' (a statement about natural numbers in~$\Eff$) is equivalent to the statement
``for any number~$n$, there \emph{\effectively} exists a prime number~$p > n$''
(a statement about natural numbers and computability in the standard topos).
(The meaning of~\emph{effectivity} will be recalled in
Section~\ref{sect:effective-topos}.)

In the syntactic view, the \effective topos merely provides a coherent way of
adding the qualification ``\effective'' to mathematical statements, for
instance turning the statement ``for any number~$n$, there exists a prime number~$p
> n$'' into ``for any number~$n$, there \emph{\effectively} exists a prime
number~$p > n$''. Similarly, a sheaf topos~$\Sh(X)$ provides a coherent way for
turning statements about real numbers and real functions into statements about
continuous~$X$-indexed families of real numbers and real functions.

The crucial point is that the translation scheme provided by any topos is sound
with respect to intuitionistic logic. Hence, regardless of our actual position
on toposes as alternate universes, working under the lens of a given topos
\emph{feels like} working in an alternate universe.

\section{The \effective topos, a universe shaped by computability}
\label{sect:effective-topos}

A basic question in computability theory is: Which computational tasks are
solvable in principle by computer programs? For instance, there is an algorithm
for computing the greatest common divisor of any pair of natural numbers, and
hence we say ``any pair of natural numbers \emph{\effectively} has a greatest
common divisor''  or ``the function~$\NN \times \NN \to \NN,\,(n,m) \mapsto
\operatorname{gcd}(n,m)$ is \emph{computable}''.

In such questions of computability, practical issues such as resource
constraints or hardware malfunctions are ignored; we employ the theoretical
notion of \emph{Turing machines}, a mathematical abstraction of the computers
of the real world.

A basic observation in computability theory is that there are computational tasks
which are not solvable even for these idealized Turing machines. The premier
example is the \emph{halting problem}: Given a Turing machine~$M$, determine
whether~$M$ terminates (comes to a stop after having carried out finitely
many computational steps) or not.

A Turing machine~$H$ which would solve this problem, that is read the
description of a Turing machine~$M$ as input and output one or zero depending
on whether~$M$ terminates or not, would be called a \emph{halting oracle}, and
a basic result is that there are no halting oracles. If we fix some \effective
enumeration of all Turing machines, then we can express the undecidability of
the halting problem also by saying that the \emph{halting function}
\[ h : \NN \longrightarrow \NN,\ n \longmapsto \begin{cases}
  1, & \text{if the~$n$-th Turing machine terminates}, \\
  0, & \text{otherwise,}
\end{cases} \]
is not computable.

The \emph{\effective topos}~$\Eff$ is a convenient home for computability
theory. A statement is true in~$\Eff$ if and only if it has a \emph{computable
witness}. For instance, a computable witness of a statement of the
form~``$\forall x\_ \exists y\_ \varphi(x,y)$'' is a Turing machine which, when
given an input~$x$, computes an output~$y$ together with a computable witness
for~$\varphi(x,y)$.

Section~\ref{sect:eff-examples} presents several examples to convey an
intuitive understanding of truth in the \effective topos; the precise
translation rules are displayed in Table~\ref{table:eff}. A precise definition
of the \effective topos requires notions of category theory which we do not
want to suppose here; it is included only for reference.

Introductory literature on the \effective topos
includes the references~\cite{hyland:effective-topos,oosten:realizability,phoa:effective,bauer:c2c}.

\begin{definition}\ \\[-1.2em]\begin{enumerate}
\item An \emph{assembly} is a set~$X$ together with a
relation~$({\realizes_X}) \subseteq \NN \times X$ such that for every element~$x
\in X$, there is a number~$n$ such that~$n \realizes_X x$.
\item A \emph{morphism of
assemblies}~$(X,{\realizes_X}) \to (Y,{\realizes_Y})$ is a map~$f : X \to Y$
which is \emph{tracked} by a Turing machine, that is for which there exists a
Turing machine~$M$ such that for any element~$x \in X$ and any number~$n$ such
that~$n \realizes x$, the computation~$M(n)$ terminates and~$M(n) \realizes
f(x)$.\end{enumerate}\end{definition}

A number~$n$ such that~$n \realizes_X x$ is called a \emph{realizer} for~$x$
and can be pictured as a concrete representation of the abstract element~$x$.
The \emph{assembly of natural numbers} is the assembly~$(\NN,{=_\NN})$ and the
\emph{assembly of functions~$\NN \to \NN$} is the
assembly~$(X,{\realizes})$ where~$X$ is the set of computable functions~$\NN
\to \NN$ and~$n \realizes f$ if and only if the~$n$-th Turing machine
computes~$f$. The category of assemblies is a regular category, but it is
missing \effective quotients. The \effective topos is obtained by a suitable
completion procedure:

\begin{definition}The \emph{\effective topos}~$\Eff$ is the ex/reg completion (as
in~\cite[Section~3.4]{menni:exact}) of the category of assemblies.\end{definition}


\subsection{Exploring the \effective topos}
\label{sect:eff-examples}
Due to its computational nature, truth in the \effective topos is quite
different from truth in the standard topos. This section explores the following
examples:

\bigskip
\begin{center}
\begin{tabular}{ll@{\qquad}l}
  \toprule
  Statement & in $\Set$ & in $\Eff$ \\
  \midrule
  Any natural number is prime or not prime. & \cmark{} (trivially) & \cmark \\
  There are infinitely many primes. & \cmark & \cmark \\
  Any function $\NN \to \NN$ is constantly zero or not. & \cmark{} (trivially) & \xmark \\
  Any function $\NN \to \NN$ is computable. & \xmark & \cmark{} (trivially) \\
  Any function $\RR \to \RR$ is continuous. & \xmark & \cmark \\
  Markov's principle holds. & \cmark{} (trivially) & \cmark \\
  Heyting arithmetic is categorical. & \xmark & \cmark \\
  \bottomrule
\end{tabular}
\end{center}
\bigskip

\newcommand{\dotparagraph}[1]{\noindent\textbf{#1}}

\begin{example}\textbf{``Any natural number is prime or not.''} Even without knowing what
a prime number is, one can safely judge this statement to be true in
the standard topos, since it is just an instance of the law of excluded middle.

By the Kripke--Joyal semantics, stating that this statement is true in the \effective topos
amounts to stating that there is a Turing machine which, when given a natural
number~$n$ as input, terminates with a correct judgment whether~$n$ is prime or
not. Such a Turing machine indeed exists -- writing such a program is often a
first exercise in programming courses. Hence the statement is also true in the
\effective topos, but for the nontrivial reason that primality can be
algorithmically tested.
\end{example}

\begin{example}\textbf{``There are infinitely many primes.''} A first-order formalization
of this statement is ``for any natural number~$n$, there is a prime
number~$p$ which is greater than~$n$'', and is known to be true in the standard
topos by any of the many proofs of this fact.

Its external meaning when interpreted in the \effective topos is that there exists
a Turing machine~$M$ which, when given a natural number~$n$ as input, terminates with a
prime number~$p > n$ as output. Such a Turing machine exists, hence the
statement is true in the \effective topos.\footnote{More precisely, the
machine~$M$ should also output the description of a Turing machine which
witnesses that~$p$ is prime and that~$p > n$. However, the statement~``$p$ is prime and~$p > n$''
is~$\neg\neg$-stable (even decidable), and for those statements witnesses are
redundant.}
\end{example}

\begin{example}\textbf{``Any function~$\NN \to \NN$ is constantly zero or not.''} Precisely,
the statement is
\[ \forall f \? \NN^\NN\_
  \bigl((\forall n \? \NN\_ f(n) = 0) \vee
  \neg
  (\forall n \? \NN\_ f(n) = 0)\bigr). \]
By the law of excluded middle, this statement is trivially true in the standard
topos.

Its meaning when interpreted in the \effective topos is that there exists a
Turing machine~$M$ which, when given the description of a Turing machine~$F$ which
computes a function~$f : \NN \to \NN$ as input, terminates with a correct
judgment of whether~$f$ is the zero function or not. Such a machine~$M$ does
not exist, hence the statement is false in the \effective topos.

Intuitively, the issue
is the following. Turing machines are able to simulate other Turing machines,
hence~$M$ could simulate~$F$ on various inputs to search the list of
function values~$f(0), f(1), \ldots$ for a nonzero number. In case that after
a certain number of steps a nonzero function value is found, the machine~$M$
can correctly output the judgment that~$f$ is not the zero function. But if the
search only turned up zero values, it cannot come to any verdict -- it cannot
rule out that a nonzero function value will show up in the as yet unexplored
part of the function.

A rigorous proof that such a machine~$M$ does not exist reduces its assumed
existence to the undecidability of the halting problem.
\end{example}

\begin{remark}Quite surprisingly, there are infinite sets~$X$ for which any flavor of constructive
mathematics, in particular the kind which is valid in any topos, verifies the
\emph{omniscience principle}
\[ \forall f\?\BB^X\_ \bigl((\exists x\?X\_ f(x) = 0) \vee (\forall x\?X\_ f(x)
= 1)\bigr), \]
where~$\BB = \{ 0, 1 \}$ is the set of booleans. This is not the case for~$X =
\NN$, but it is for instance the case for the one-point compactification~$X =
\NN_\infty$ of the naturals. This phenomenon has been thoroughly explored by
Escardó~\cite{escardo:omniscience}.\end{remark}

\begin{example}\textbf{``Any function~$\NN \to \NN$ is computable.''}
The preceding examples give the impression that what is true in the
\effective topos is solely a subset of what is true in the standard topos. The
example of this subsection, the so-called \emph{formal Church--Turing thesis},
shows that the relation between the two toposes is more nuanced.

As recalled above, in the standard topos there are functions~$\NN \to \NN$ which are not computable by a Turing
machine. Cardinality arguments
even show that most functions~$\NN \to \NN$ are not computable: There
are~$\aleph_0^{\aleph_0} = 2^{\aleph_0}$ functions~$\NN \to \NN$, but
only~$\aleph_0$ Turing machines and hence only~$\aleph_0$ functions which are
computable by a Turing machine.

In contrast, in the \effective topos, any function~$\NN \to \NN$ is computable
by a Turing machine. The external meaning of this internal statement is that
there exists a Turing machine~$M$ which, when given a description of a Turing
machine~$F$ computing a function~$f : \NN \to \NN$, outputs a description of a
Turing machine computing~$f$. It is trivial to program such a machine~$M$: the
machine~$M$ simply has to echo its input back to the user.

To avert a paradox, we should point out where the usual proof of the
existence of noncomputable functions theory employs nonconstructive reasoning, for if
the proof would only use intuitionistic reasoning, it would also hold internally to the
\effective topos, in contradiction to the fact that in the \effective topos
all functions~$\NN \to \NN$ are computable.

The usual proof sets up the halting function~$h : \NN \to \NN$, defined using
the case distinction
\[ h : n \mapsto \begin{cases}
  1, & \text{if the $n$-th Turing machine terminates}, \\
  0, & \text{if the $n$-th Turing machine does not terminate},
\end{cases} \]
and proceeds to show that~$h$ is not computable. However, in the \effective
topos, this definition does not give rise to a total function from~$\NN$
to~$\NN$. The actual domain is the subset~$M$ of those natural numbers~$n$
for which the~$n$-th Turing machine
terminates or does not terminate. This condition is trivial only assuming the
law of excluded middle; intuitionistically, this condition is nontrivial and
cuts out a nontrivial subset of~$\NN$.

Subobjects in the \effective topos are more than mere subsets; to give an
element of~$M$ in the \effective topos, we need not only give a natural
number~$n$ such that the~$n$-th Turing machine terminates or does not
terminate, but also a computational witness of either case. For any particular
numeral~$n_0 \in \NN$, there is such a witness (appealing to the law of
excluded middle on the meta level), and hence the statement~``$n_0
\in M$'' holds in the \effective topos. However, there is no program which
could compute such witnesses for any number~$n$, hence the statement~``$\forall
n\?\NN\_ n \in M$'' is not true in~$\Eff$ and hence the \effective topos does not believe~$M$
and~$\NN$ to be the same.
\end{example}

\begin{example}\textbf{``Any function~$\RR \to \RR$ is continuous.''}\label{sect:eff-continuous}
In the standard topos, this statement is plainly false, with the sign and Heaviside step functions
being prominent counterexamples. In the \effective topos, this statement is
true and independently due to~\cite{kreisel-lacombe-shoenfield:cont} and~\cite{ceitin:cont}. A rigorous proof is not entirely
straightforward (a textbook reference
is~\cite[Theorem~9.2.1]{longley-normann:higher-order-computability}), but an intuitive
explanation is as follows.

What the \effective topos believes to be a real number is, from the external
point of view, a Turing machine~$X$ which outputs, when called with a natural
number~$n$ as input, a rational approximation~$X(n)$. These approximations are
required to be \emph{consistent} in the sense that~$|X(n) - X(m)|
\leq 2^{-n} + 2^{-m}$. Intuitively, such a machine~$X$ denotes the real
number~$\lim_{n \to \infty} X(n)$, and each approximation~$X(n)$ must be
within~$2^{-n}$ of the limit.

A function~$f : \RR \to \RR$ in the \effective topos is therefore given by a
Turing machine~$M$ which, when given the description of such a Turing machine~$X$ as
input, outputs the description of a similar such Turing machine~$Y$.
To compute a rational approximation~$Y(n)$, the machine~$Y$ may simulate~$X$
and can therefore determine arbitrarily many rational approximations~$X(m)$.
However, within a finite amount of time, the machine~$Y$ can only learn finitely many
such approximations. Hence a function such as the sign function, for which
even rough rational approximations of~$\operatorname{sgn}(x)$ require infinite
precision in the input~$x$, does not exist in the \effective topos.
\end{example}

\begin{example}\textbf{``Markov's principle holds.''} Markov's principle is the
following statement:
\begin{equation}\label{eq:markov}\tag{MP}
  \forall f \? \NN^\NN\_ \bigl((\neg\neg\exists n\?\NN\_ f(n) = 0)
  \Longrightarrow \exists n\?\NN\_ f(n) = 0\bigr).
\end{equation}
It is an instance of the law of double negation elimination and hence trivially
true in the standard topos, at least if we subscribe to classical logic on the
meta level. A useful consequence of Markov's principle is that Turing machines
which do not run forever (that is, which do \emph{not not} terminate) actually
terminate; this follows by applying Markov's principle to the function~$f : \NN
\to \NN$ where~$f(n)$ is zero or one depending on whether a given Turing
machine has terminated within~$n$ computational steps or not.

The \effective topos inherits Markov's principle from the meta level:
The statement~``$\Eff \models \text{\eqref{eq:markov}}$'' means that there is a
Turing machine~$M$ which, when given the description of a Turing machine~$F$
computing a function~$f : \NN \to \NN$ as input, outputs the description of a Turing
machine~$S_F$ which, when given a witness of~``$\neg\neg \exists n\?\NN\_ f(n)
= 0$'', outputs a witness of~``$\exists n\?\NN\_ f(n) = 0$'' (up to trivial
conversions, this is a number~$n$ such that~$f(n) = 0$).

By the translation rules listed in Table~\ref{table:eff}, a number~$e$
realizes~``$\neg\neg \exists n\?\NN\_ f(n) = 0$'' if and only if it is
\emph{not not} the case that there is some number~$e'$ such that~$e'$
realizes~``$\exists n\?\NN\_ f(n) = 0$''. Hence, if ``$\exists n\?\NN\_ f(n) =
0$'' is realized at all, then any number is a witness of~``$\neg\neg
\exists n\?\NN\_ f(n) = 0$''.

As a consequence, the input given to the machine~$S_F$ is entirely uninformative
and~$S_F$ cannot make direct computational use of it. But its existence ensures that an
\emph{unbounded search} will not fail (and hence succeed, by an appeal to Markov's
principle on the meta level): The machine~$S_F$ can simulate~$F$ to
compute the values~$f(0), f(1), f(2), \ldots$ in turn, and stop
with output~$n$ as soon as it determines that some function value~$f(n)$ is zero.
\end{example}

\begin{example}\textbf{``Heyting arithmetic is categorical.''} In addition to the standard
model~$\NN$, the standard topos contains uncountably many nonstandard models of
Peano arithmetic (at least if we assume a classical meta level). By a theorem
of~\cite{berg-oosten:arithmetic}, the situation
is quite different in the \effective topos:
\begin{enumerate}
\item Heyting arithmetic, the intuitionistic cousin of Peano arithmetic, is
categorical in the sense that it has
exactly one model up to isomorphism, namely~$\NN$.
\item In fact, even the finitely axiomatizable subsystem of Heyting arithmetic
where the induction scheme is restricted to~$\Sigma_1$-formulas has exactly one
model up to isomorphism, again~$\NN$. As a consequence, Heyting arithmetic is
finitely axiomatizable.
\item Peano arithmetic is ``quasi-inconsistent'' in that it does not have any
models, for any model of Peano arithmetic would also be a model of Heyting
arithmetic, but the only model of Heyting arithmetic is~$\NN$ and~$\NN$ does
not validate the theorem ``any Turing machine terminates or does not
terminate'' of Peano arithmetic.
\end{enumerate}
As a consequence, Gödel's completeness theorem fails in the \effective topos:
In the \effective topos, Peano arithmetic is consistent (because it is
equiconsistent to Heyting arithmetic, which has a model) but does not have a
model.

Statement~(1) is reminiscent of the fact due to~\cite{tennenbaum:models} that no nonstandard model of
Peano arithmetic in the standard topos is computable.
\end{example}

\begin{table}
  \centering
  \begin{framed}\begin{tabbing}
    $e \models (\forall f\?\NN^\NN\_ \varphi(n))$ \= \kill
    $\Eff \models \varphi$ \> iff there is a natural number~$e$ such that~$e
    \realizes \varphi$. \\\\
    \begin{minipage}{1.04\textwidth}
    A number~$e$ such that~$e \realizes \varphi$ is called a
    \emph{realizer} for~$\varphi$. It is the precise version of what is called
    \emph{computational witness} in the main text.
    In the following, we write~``$e \cdot n \downarrow$'' to mean that
    the~$e$-th Turing machine terminates on input~$n$, and in this case denote
    the result by~``$e \cdot n$''. No separate clause for negation is listed,
    as~``$\neg\varphi$'' is an abbreviation for~``$(\varphi \Rightarrow
    \bot)$''.\end{minipage} \\\\
    $e \realizes s = t$ \> iff $s = t$. \\
    $e \realizes \top$ \> iff $1 = 1$. \\
    $e \realizes \bot$ \> iff $1 = 0$. \\
    $e \realizes (\varphi \wedge \psi)$ \> iff~$e \cdot 0 \downarrow$ and~$e
    \cdot 1 \downarrow$ and $e\cdot0 \realizes \varphi$ and~$e\cdot1 \realizes \psi$. \\
    $e \realizes (\varphi \vee \psi)$ \> iff~$e \cdot 0 \downarrow$ and~$e
    \cdot 1 \downarrow$ and \\ \> \qquad if~$e\cdot0 = 0$ then~$e\cdot1 \realizes
    \varphi$, and \\ \> \qquad if~$e\cdot0 \neq 0$ then~$e\cdot1 \realizes \psi$. \\
    $e \realizes (\varphi \Rightarrow \psi)$ \> iff for any number~$r \in \NN$
    such that~$r \realizes \varphi$, $e \cdot r \downarrow$ and~$e \cdot r \realizes \psi$. \\
    $e \realizes (\forall n\?\NN\_ \varphi(n))$ \> iff for any number~$n_0
    \in \NN$, $e \cdot n_0 \downarrow$ and~$e \cdot n_0 \realizes \varphi(n_0)$. \\
    $e \realizes (\exists n\?\NN\_ \varphi(n))$ \> iff~$e\cdot0 \downarrow$ and~$e\cdot1 \downarrow$
    and~$e\cdot1 \realizes \varphi(e\cdot0)$. \\
    $e \realizes (\forall f\?\NN^\NN\_ \varphi(f))$ \> iff for any function~$f_0
    : \NN \to \NN$ and any number~$r_0$ such that \\ \> \qquad $f_0$ is computed by the~$r_0$-th
    Turing machine, \\ \> \qquad
    $e \cdot r_0 \downarrow$ and~$e \cdot r_0 \realizes \varphi(f_0)$. \\
    $e \realizes (\exists f\?\NN^\NN\_ \varphi(f))$ \> iff~$e \cdot 0 \downarrow$
    and~$e \cdot 1 \downarrow$ and
    the $(e \cdot 0)$-th Turing machine \\ \> \qquad computes a function~$f_0 : \NN \to \NN$
    and $e \cdot 1 \realizes \varphi(f_0)$.
  \end{tabbing}\end{framed}
  \bigskip

  \caption{\label{table:eff}A (fragment of) the translation
  rules defining the meaning of statements internal to the \effective topos.}
\end{table}

\subsection{Variants of the \effective topos} The \effective topos belongs to a
wider class of \emph{realizability toposes}. These can be obtained by repeating
the construction of the \effective topos with any other reasonable model of
computation in place of Turing machines. The resulting toposes will in general
not be equivalent and reflect higher-order properties of the employed models.
Two of these further toposes are of special philosophical interest.

\bigskip
\paragraph{Hypercomputation}
Firstly, in place of ordinary Turing machines, one can employ the
\emph{infinite-time Turing machines} pioneered by~\cite{hamkins-lewis:ittm}. These machines model \emph{hypercomputation}
in that they can run for ``longer than infinity''; more precisely, their
computational steps are indexed by the ordinal numbers instead of the natural
numbers. For instance, an infinite-time Turing machine can trivially decide the
twin prime conjecture, by simply walking along the natural number line and
recording any twin primes it finds. Then, on day~$\omega$, it can observe
whether it has found infinitely many twins or not.

In the realizability topos constructed using infinite-time Turing machines, the full
law of excluded middle still fails, but some instances which are wrong in the
\effective topos do hold in this topos. For example, the instance ``any
function~$\NN \to \NN$ is the zero function or not'' does: Its external meaning
is that there is an infinite-time Turing machine~$M$ which, when given the
description of an infinite-time Turing machine~$F$ computing a function~$f :
\NN \to \NN$ as input, terminates (at some ordinal time step) with a correct
judgment of whether~$f$ is the zero function or not. Such a machine~$M$ indeed
exists: It simply has to simulate~$F$ on all inputs~$0,1,\ldots$ and
check whether one of the resulting function values is not zero. This search
will require a transfinite amount of time (not least because simulating~$F$ on
just one input might require a transfinite amount of time), but infinite-time
Turing machines are capable of carrying out this procedure.

The realizability topos given by infinite-time Turing machines provides an intriguing environment challenging many
mathematical intuitions shaped by classical logic. For instance, while from the
point of view of this topos the reals are still uncountable in the sense that
there is no surjection~$\NN \to \RR$, there is an injection~$\RR \to
\NN$~\cite[Section~4]{bauer:injection}.\footnote{What the realizability topos
given by infinite-time Turing machines believes to be a real number is, from
the external point of view, an infinite-time Turing machine~$X$ which outputs,
when called with a natural number~$n$ as input, a rational
approximation~$X(n)$. As with the original \effective topos, these
approximations have to be consistent in the sense that~$|X(n) - X(m)| \leq
2^{-n} + 2^{-m}$, and two such machines~$X, X'$ represent the same real
iff~$|X(n) - X'(m)| \leq 2^{-n} + 2^{-m}$ for all natural numbers~$n,m$.

A map~$\RR \to \NN$ in this topos is hence given by an infinite-time
Turing machine~$M$ which, when given the description of such an infinite-time
Turing machine~$X$ as input, outputs a certain natural number~$M(X)$. If~$X$
and~$X'$ represent the same real, then~$M(X)$ has to coincide with~$M(X')$.
This map~$\RR \to \NN$ is injective iff conversely~$M(X) = M(X')$ implies
that~$X$ and~$X'$ represent the same real.

We can program such a machine~$M$ as follows: Read the description of an
infinite-time Turing machine~$X$ representing a real number as input. Then
simulate, in a dovetailing fashion, all infinite-time Turing machines and
compare their outputs with the outputs of~$X$. As soon as a machine~$X'$ is
found which happens to terminate on all inputs in such a way that~$|X(n) -
X'(m)| \leq 2^{-n} + 2^{-m}$ for all natural numbers~$n,m$, output the number
of this machine (in the chosen enumeration of all infinite-time Turing
machines) and halt.

The number~$M(X)$ computed by~$M$ depends on the input/output behaviour of~$X$,
the chosen ordering of infinite-time Turing machines, and on details of the
interleaving simulation and the comparison procedure -- but it does not depend on
the implementation of~$X$ or on its specific choice of rational
approximations~$X(n)$. The search terminates since there is at least one
infinite-time Turing machine which represents the same real number as~$X$ does,
namely~$X$ itself.}

\bigskip
\paragraph{Machines of the physical world} A second variant of the
\effective topos is obtained by using machines of the physical world
instead of abstract Turing machines. In doing so, we of
course leave the realm of mathematics, as real-world machines are not objects
of mathematical study. But it is still interesting to see which commitments
about the nature of the physical world imply which internal statements of the
resulting topos.

For instance, \cite{bauer:int-mathematics} showed that inside this topos any function~$\RR \to
\RR$ is continuous if, in the physical world, only finitely many computational
steps can be carried out in finite time and if it is possible to form
tamper-free private communication channels.

\section{Toposes of sheaves, a convenient home for local truth}
\label{sect:sheaf-toposes}

Associated to any topological space~$X$ (such as Euclidean space), there is the
\emph{topos of sheaves over~$X$}, $\Sh(X)$. To a first approximation, a
statement is true in~$\Sh(X)$ if and only if it ``holds locally on~$X$'';
what~$\Sh(X)$ believes to be a set is a ``continuous family of sets, one set
for each point of~$X$''. The precise rules of the Kripke--Joyal semantics
of~$\Sh(X)$ are listed in Table~\ref{table:sheaf}.

Just as the \effective topos provides a coherent setting for studying
computability using a naive element-based language, the sheaf topos~$\Sh(X)$ provides a
coherent setting for studying continuous~$X$-indexed families of objects (sets,
numbers, functions) as if they were single objects.

Sheaf toposes take up a special place in the history of topos theory: If the
base~$X$ is allowed to be a \emph{site} instead of a topological space, the
resulting toposes constitute the large class of Grothendieck toposes, the
original notion of toposes. Categorically, the passage from topological spaces
to sites is rather small, but the resulting increase in flexibility is
substantial and fundamental to modern algebraic geometry.

\subsection{A geometric interpretation of double negation}
In intuitionistic logic, the double negation~$\neg\neg\varphi$ of a
statement~$\varphi$ is a slight weakening of~$\varphi$; while~$(\varphi
\Rightarrow \neg\neg\varphi)$ is an intuitionistic tautology, the converse can
only be shown for some specific statements. The internal language of~$\Sh(X)$
gives geometric meaning to this logical peculiarity.

Namely, it is an instructive exercise that~$\Sh(X) \models \neg\neg\varphi$ is equivalent to the
existence of a \emph{dense open}~$U$ of~$X$ such that~$U \models \varphi$.
If~$\Sh(X) \models \varphi$, that is if~$X \models \varphi$, then there
obviously exists such a dense open, namely~$X$ itself; however the converse
usually fails.

The only case that the law of excluded middle does hold internally to~$\Sh(X)$
is when the only dense open of~$X$ is~$X$ itself; assuming classical logic in
the metatheory, this holds if and only if every open is also closed. This is
essentially only satisfied if~$X$ is discrete.

An important special case is when~$X$ is the one-point space. In this
case~$\Sh(X)$ is equivalent (as categories and hence toposes) to the standard
topos. To the extent that mathematics within~$\Sh(X)$ can be described as ``mathematics
over~$X$'', this observation justifies the slogan that ``ordinary mathematics
is mathematics over the point''.

\subsection{Reifying continuous families of real numbers as single real numbers}
As detailed in Section~\ref{sect:eff-continuous}, what the \effective topos believes to be a real number is
actually a Turing machine computing arbitrarily-good consistent rational
approximations. A similarly drastic shift in meaning, though in an orthogonal
direction, occurs with~$\Sh(X)$. What~$\Sh(X)$ believes to be a (Dedekind) real
number~$a$ is actually a continuous family of real numbers on~$X$, that is, a
continuous function~$a : X \to \RR$~\cite[Corollary~D4.7.5]{johnstone:elephant}.

Such a function is everywhere positive on~$X$ if and only if, from the internal point of
view of~$\Sh(X)$, the number~$a$ is positive; it is everywhere zero if and only
if, internally, the number~$a$ is zero; and it is everywhere negative if and
only if, internally, the number~$a$ is negative.

\begin{figure}
  \centering
  \includegraphics[height=3cm]{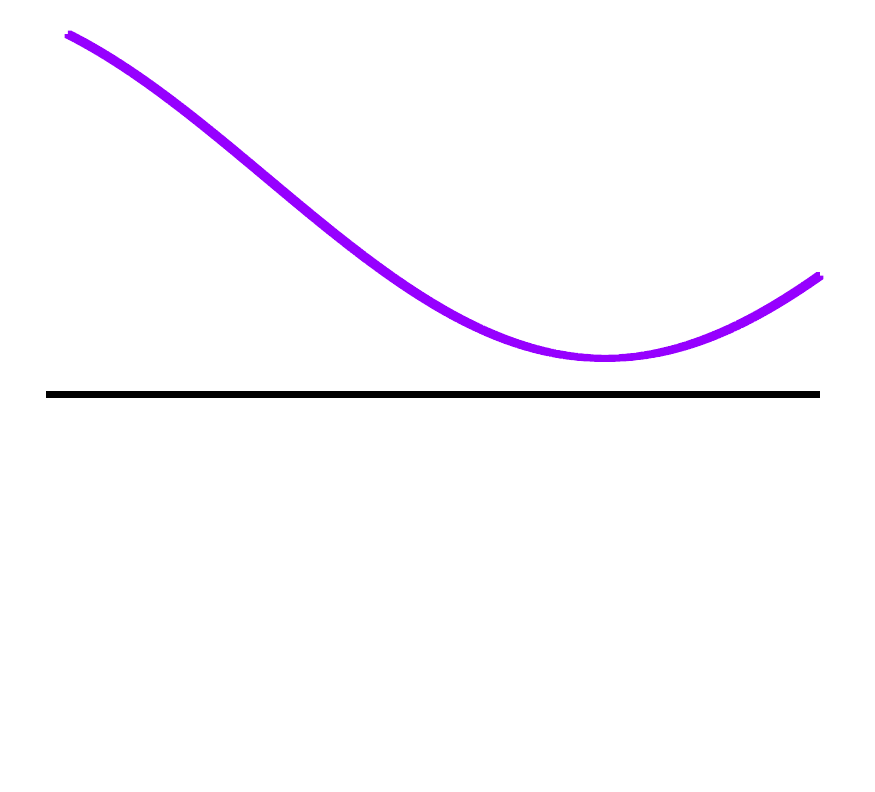}
  \includegraphics[height=3cm]{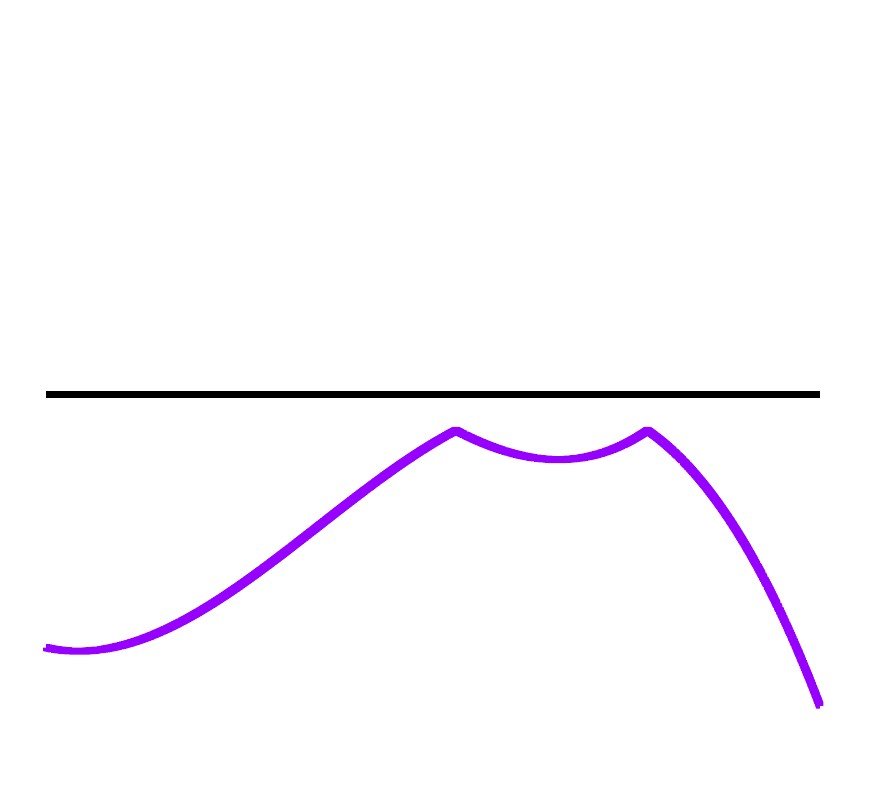}
  \includegraphics[height=3cm]{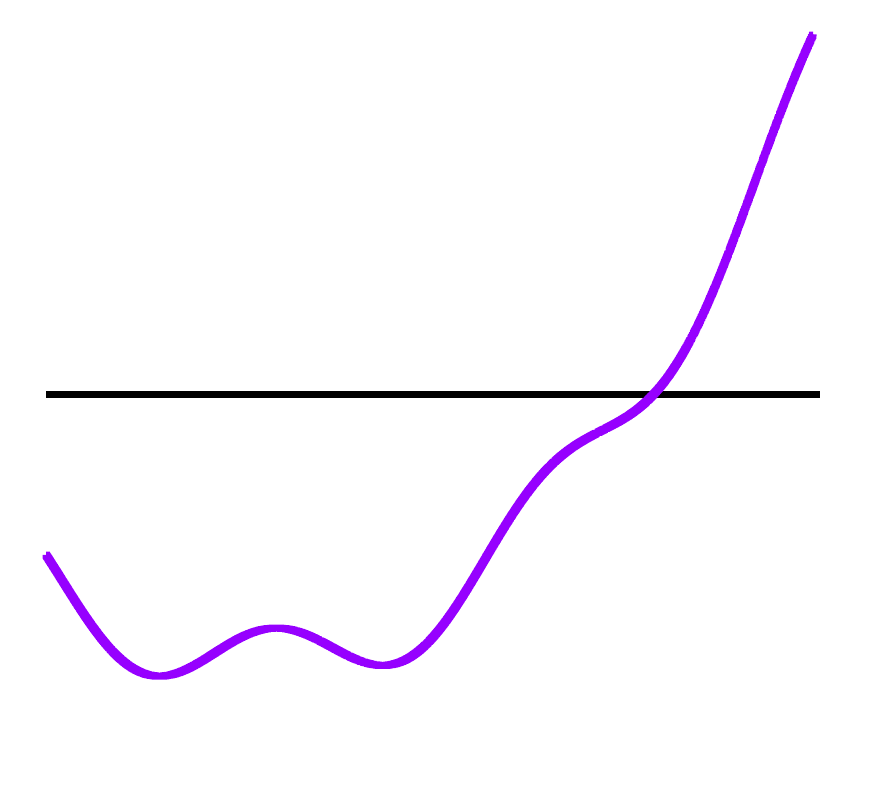}
  \caption{\label{fig:trichotomy}Three examples of what the topos~$\Sh(X)$
  believes to be a single real number, where the base space~$X$ is the
  unit interval. (a) A positive real number. (b) A negative real number. (c)~A~number
  which is neither negative nor zero nor positive. Externally speaking,
  there is no covering of the unit interval by opens on which the
  depicted function~$a$ is either everywhere negative, everywhere zero or everywhere
  positive.}
\end{figure}

The law of trichotomy, stating that any real number is either negative, zero or
positive, generally fails in~$\Sh(X)$. By the Kripke--Joyal semantics, the external
meaning of the internal statement~``$\forall a:\RR\_ a<0 \vee a=0 \vee a>0$'' is that for any continuous function~$a : U
\to \RR$ defined on any open~$U$ of~$X$, there is an open covering~$U =
\bigcup_i U_i$ such that on each member~$U_i$ of this covering, the function~$a$ is
either everywhere negative on~$U_i$, everywhere zero on~$U_i$ or everywhere
positive on~$U_i$. But this statement is, for most base spaces~$X$, false.
Figure~\ref{fig:trichotomy}(c) shows a counterexample.

The weaker statement ``for any real number~$a$ it is \notnot the case that~$a < 0$
or~$a = 0$ or~$a > 0$'' does hold in~$\Sh(X)$, for this statement is a
theorem of intuitionistic calculus. Its meaning is that there exists a
dense open~$U$ such that~$U$ can be covered by opens on which~$a$ is either
everywhere negative, everywhere zero or everywhere positive. In the example
given in Figure~\ref{fig:trichotomy}(c), this open~$U$ could be taken as~$X$
with the unique zero of~$a$ removed.

\begin{table}
  \centering
  \begin{framed}\begin{tabbing}
    $U \models (\forall x\?\RR\_ \varphi(x))$ \= \kill
    $\Sh(X) \models \varphi$ \> iff $X \models \varphi$. \\\\
    $U \models a = b$ \> iff~$a = b$ on~$U$. \\
    $U \models \top$ \> is true for any open~$U$. \\
    $U \models \bot$ \> iff~$U$ is the empty open. \\
    $U \models (\varphi \wedge \psi)$ \> iff~$U \models \varphi$ and~$U \models \psi$. \\
    $U \models (\varphi \vee \psi)$ \> iff there is an open covering~$U =
    \bigcup_i U_i$ such that, \\ \> \qquad for each index~$i$, $U_i \models \varphi$
    or $U_i \models \psi$. \\
    $U \models (\varphi \Rightarrow \psi)$ \> iff for any open~$V \subseteq U$,
    $V \models \varphi$ implies~$V \models \psi$. \\
    $U \models (\forall a\?\RR\_ \varphi(a))$ \> iff for any open~$V
    \subseteq U$ and \\\>\qquad any continuous function~$a_0 : V \to \RR$, $V \models
    \varphi(a_0)$. \\
    $U \models (\exists a\?\RR\_ \varphi(a))$ \> iff there is an open
    covering~$U = \bigcup_i U_i$ such that, \\ \> \qquad for each index~$i$,
    there exists a \\\>\qquad continuous function~$a_0 : U_i \to \RR$ with~$U_i \models
    \varphi(a_0)$.
  \end{tabbing}\end{framed}
  \bigskip

  \caption{\label{table:sheaf}A (fragment of) the translation rules defining
  the meaning of statements internal to~$\Sh(X)$, the topos of sheaves over a
  topological space~$X$.}
\end{table}

\begin{figure}
  \centering
  \marginnote{
    \begin{ocg}[printocg=never]{layername}{layerid1}{on}
      \vspace*{-14em}
      \scalebox{0.2}{
        \animategraphics[loop,autoplay]{60}{zeros-in-families-frames/}{0000}{0188}
      }
      \begin{minipage}{3cm}\scriptsize\noindent
      In the online version of this book, a video of the
      continuous family is shown here.\end{minipage}
    \end{ocg}
  }
  \includegraphics[height=3cm]{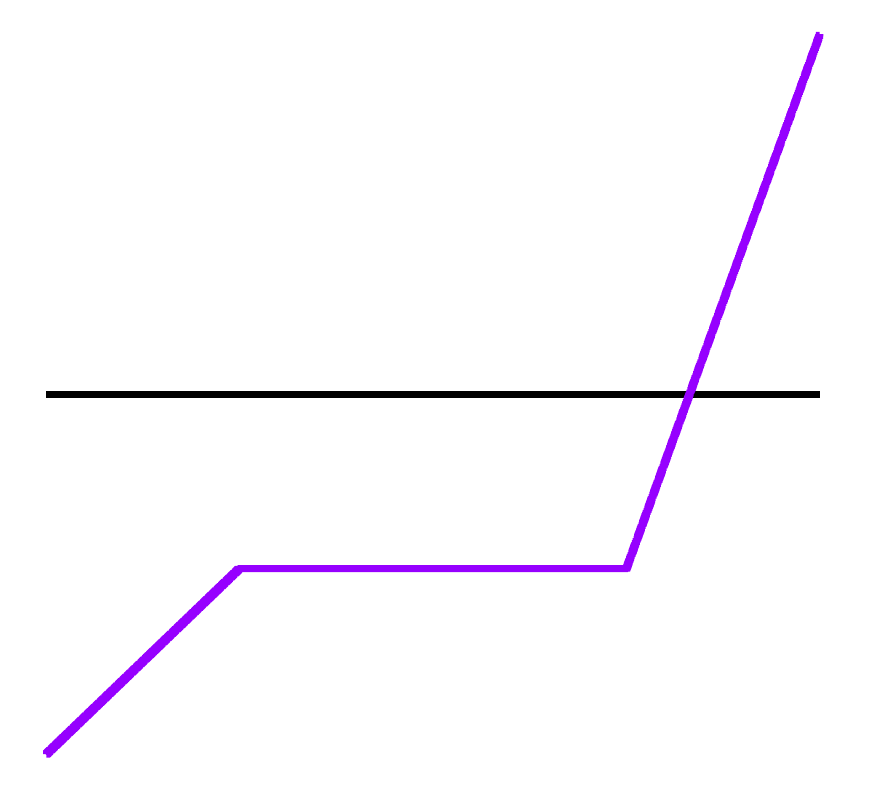}
  \includegraphics[height=3cm]{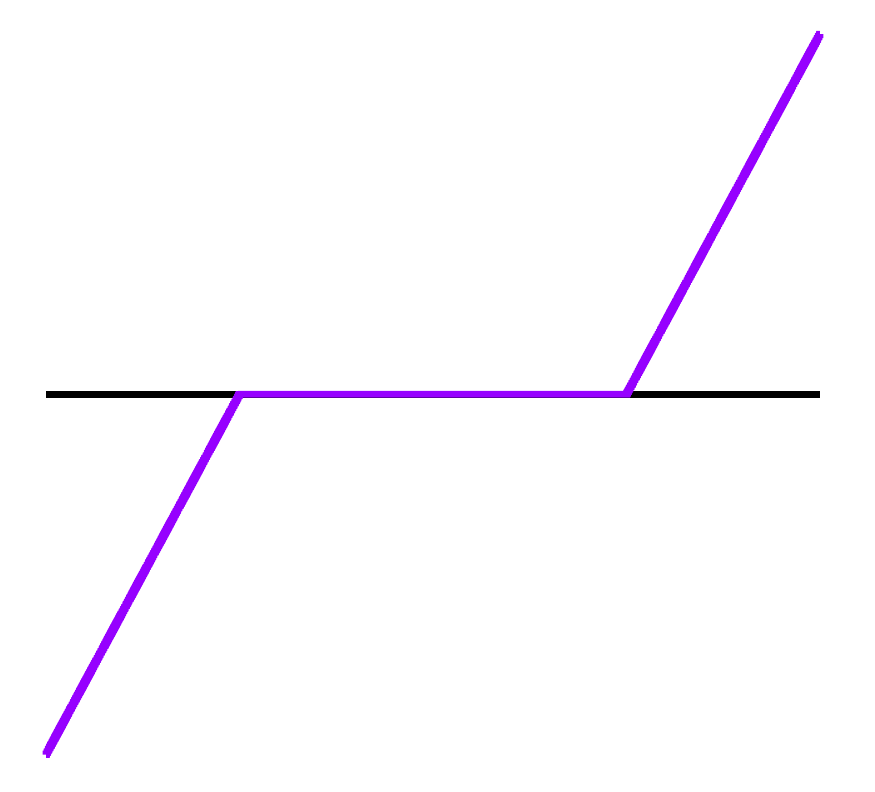}
  \includegraphics[height=3cm]{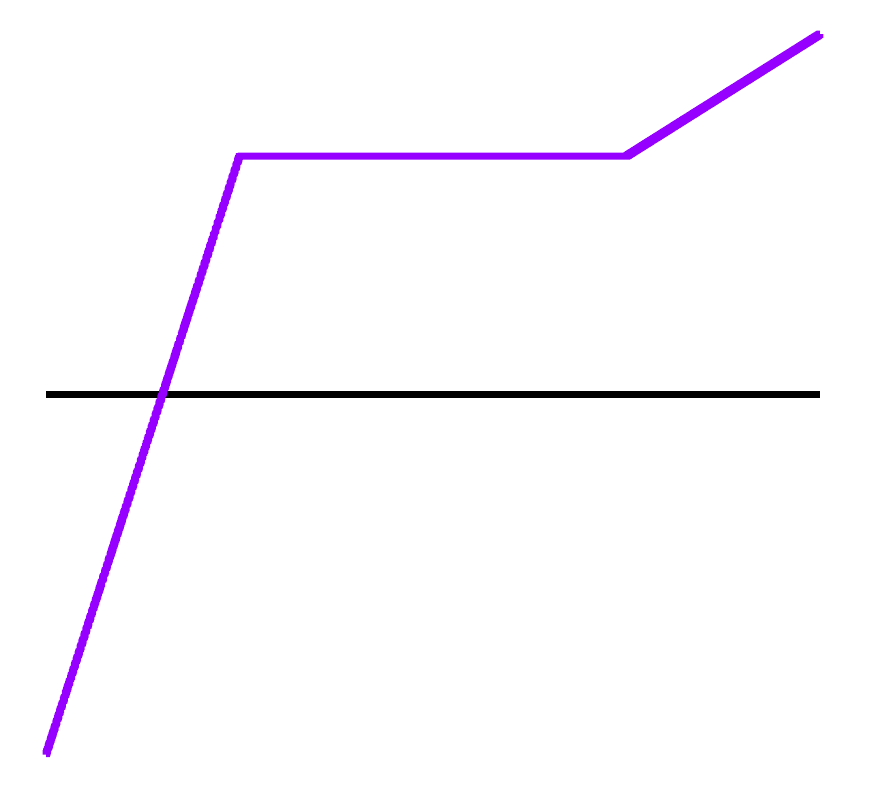}

  \caption{\label{fig:ivt}Three members~$f_{x_0}, f_{x_1}, f_{x_2}$ of a continuous
  family~$(f_x)_{x \in X}$ of continuous functions~$f_x : \RR \to \RR$. The
  parameter space is~$X = [0,1]$ (not shown). The functions~$f_x$ are
  obtained by moving the horizontal plateau up or
  down. The leftmost depicted member~$f_{x_0}$ has a unique
  zero, and there is an open neighborhood~$U$ of~$x_0$ on which zeros of the
  functions~$f_x$, $x \in U$ can be picked continuously. The same is true
  for~$x_2$ (right figure). However, there is no
  such neighborhood of that particular parameter value~$x_1$ for which the
  horizontal plateau lies on the $x$-axis (middle figure).}
\end{figure}

\subsection{Reifying continuous families of real functions as single real functions}
Let~$(f_x)_{x \in X}$ be a continuous family of
continuous real-valued functions; that is, not only should each of the individual
functions~$f_x : \RR \to \RR$ be continuous, but the joint map~$\RR
\times X \to \RR, (a,x) \mapsto f_x(a)$ should be continuous. (This stronger
condition implies continuity of the individual functions.)
From the point of view of~$\Sh(X)$, this family looks like a single continuous
function~$f : \RR \to \RR$.

The internal statement~``$f(-1) < 0$'' means that~$f_x(-1) < 0$ for all~$x \in
X$, and similarly so for being positive. More generally, if~$a$ and~$b$ are
continuous functions~$X \to \RR$ (hence real numbers from the internal point of
view), the internal statement~``$f(a) < b$'' means that~$f_x(a(x)) < b(x)$ for
all~$x \in X$.

The internal statement~``$f$ possesses a zero'', that is ``there exists a
number~$a$ such that~$f(a) = 0$'', means that all the functions~$f_x$ each
possess a zero and that moreover, these zeros can locally be picked in a
continuous fashion. More precisely, this statement means that there is an open
covering~$X = \bigcup_i U_i$ such that, for each index~$i$, there is a continuous
function~$a : U_i \to \RR$ such that~$f_x(a(x)) = 0$ for all~$x \in U_i$. (On
overlaps~$U_i \cap U_j$, the zero-picking functions~$a$ need not agree.)

\begin{example}
From these observations we can deduce that the intermediate value theorem of
undergraduate calculus does in general not hold in~$\Sh(X)$ and hence does not allow
for an intuitionistic proof. The intermediate value theorem states: ``If~$f : \RR \to \RR$
is a continuous function such that~$f(-1) < 0$ and~$f(1) > 0$, there exists a
number~$a$ such that~$f(a) = 0$.'' The external meaning of this statement is
that for any continuous family~$(f_x)_x$ of continuous functions with~$f_x(-1) <
0$ and~$f_x(1) > 0$ for all~$x \in X$, it is locally possible to pick zeros of
the family in a continuous fashion. Figure~\ref{fig:ivt} shows a counterexample
to this claim.

In contrast, the intermediate value theorem for (strictly) monotone functions does have an
intuitionistic proof and hence applies in the internal universe of~$\Sh(X)$.
Thus for any continuous family~$(f_x)_x$ of continuous monotone functions
with~$f_x(-1) < 0$ and~$f_x(1) > 0$ for all~$x \in X$, it is locally possible
to pick zeros of the family in a continuous fashion.\footnote{While the
Kripke--Joyal translation of~``$\exists$'' is by definition \emph{local
existence}, one can show that the Kripke--Joyal translation of~``$\exists!$''
is \emph{unique existence on all opens}, in particular \emph{unique global
existence}. Because the conclusion of the intermediate value theorem for
monotone functions can be strengthened from ``has a zero'' to ``has a unique
zero'', this observation shows that the zeros can even globally be picked in a
continuous fashion.}
\end{example}

\begin{example}The fundamental theorem of algebra generally fails in~$\Sh(X)$, even
for quadratic polynomials. What~$\Sh(X)$ believes to be a (Dedekind)
complex number is externally a continuous function~$X \to \CC$. Let~$X$ be the
complex plane. Then the identity function~$\mathrm{id}_X$ is a single complex number from the
internal point of view of~$\Sh(X)$. The fundamental theorem of algebra would
predict~``$\Sh(X) \models \exists a\?\CC\_ a^2 - \mathrm{id}_X = 0$'',
hence that there is an open covering~$X = \bigcup_i U_i$ such that on each
open~$U_i$, there is a continuous function~$a : U_i \to \CC$ such that~$a(z)^2
- z = 0$ for all~$z \in U_i$. However, it is a basic fact of complex analysis
that such a function does not exist if~$0 \in U_i$.
\end{example}

\begin{example}The standard proof of Banach's fixed point theorem employs only
intuitionistic reasoning, hence applies internally to~$\Sh(X)$. Interpreting the
internal Banach fixed point theorem by the Kripke--Joyal translation rules
yields the statement that fixed points of continuous families of contractions
depend continuously on parameters.
\end{example}

\section{Toposes adapted to synthetic differential geometry}
\label{sect:smooth}

The idea of \emph{infinitesimal numbers} -- numbers which can be pictured as lying between~$-\frac{1}{n}$
and~$\frac{1}{n}$ for any natural number~$n$ (though this intuition will not
serve as their formal definition in this text) -- has a long and rich history. They are
not part of today's standard setup of the reals, but they are still intriguing
as a calculational tool and as a device to bring mathematical intuition and
mathematical formalism closer together.

For instance, employing numbers~$\varepsilon$ such that~$\varepsilon^2 = 0$, we can
compute derivatives blithely as follows, without requiring the notion of
limits:
\begin{equation}\label{eq:derivative}\tag{$\star$}
  \begin{aligned}
    (x + \varepsilon)^2 - x^2 &= x^2 + 2x\varepsilon + \varepsilon^2 - x^2 = 2x\varepsilon \\
    (x + \varepsilon)^3 - x^3 &= x^3 + 3x^2\varepsilon + 3x\varepsilon^2 + \varepsilon^3 - x^3 = 3x^2\varepsilon
  \end{aligned}
\end{equation}
In each case, the derivative is visible as the coefficient of~$\varepsilon$ in
the result. A further example is from geometry: Having a nontrivial
set~$\Delta$ of infinitesimal numbers available allows us to define a
\emph{tangent vector} to a manifold~$M$ to be a map~$\gamma : \Delta \to M$. This
definition precisely captures the intuition that a tangent vector is an
infinitesimal curve.

\subsection{Hyperreal numbers} There are several ways of introducing
infinitesimals into rigorous mathematics. One is Robinson's \emph{nonstandard
analysis}, where we enlarge the field~$\RR$ of real numbers to a
field~$^\star\RR$ of \emph{hyperreal numbers} by means of a non-principal
ultrafilter.

The hyperreals contain an isomorphic copy of the ordinary reals as the
so-called \emph{standard elements}, and they also contain infinitesimal numbers
and their inverses, transfinite numbers. Additionally, they support a powerful \emph{transfer
principle}: Any statement which does not refer to standardness is true for the
hyperreals if and only if it is true for the ordinary reals.

In the ``if'' direction, the transfer principle is useful for importing
knowledge about the ordinary reals into the hyperreal realm. For instance,
addition of hyperreals is commutative because addition of reals is.
By the ``only if'' direction, a theorem established for the hyperreals also
holds for the ordinary reals. In this way, the infinitesimal numbers of
nonstandard analysis can be viewed as a convenient fiction, generating a
conservative extension of the usual setup of mathematics.

There is a growing body of research in mathematics which employs hyperreal
numbers in this sense. To exemplarily cite just one example, a recent
application of nonstandard analysis in symplectic geometry is due to~\cite{fabert:non-squeezing,fabert:floer}, who verified an
infinite-dimensional analogue of the Arnold conjecture.

However, the realization of the fiction of infinitesimal numbers in nonstandard analysis crucially rests on a non-principal
ultrafilter, whose existence requires principles which go beyond the means of
Zermelo--Fraenkel set theory~\textsc{zf}.\footnote{A hyperreal number is
represented by an infinite sequence~$(x_0,x_1,x_2,\ldots)$ of ordinary real
numbers. For instance, the sequence~$(1,1,1,\ldots)$ represents the hyperreal
version of the number~$1$, the sequence~$(1,\frac{1}{2},\frac{1}{3},\ldots)$
represents an infinitesimal number and its inverse~$(1,2,3,\ldots)$ represents
a transfinite number.
The sequence~$(1,1,1,\ldots)$ is deemed positive, and so is~$(-1,1,1,1,\ldots)$,
which differs from the former only in finitely many places. But
should~$(1,-1,1,-1,\ldots)$ be deemed positive or negative? Whatever the
answer, our decision has consequences for other sequences. For
instance~$(-1,1,-1,1,\ldots)$ should be assigned the opposite sign
and~$(\tan(1),\tan(-1),\tan(1),\tan(-1),\ldots)$ the same.
A non-principal ultrafilter is a set-theoretic gadget which fixes all such
decisions once and for all in a coherent manner. Having such an ultrafilter
available, a sequence~$(x_0,x_1,x_2,\ldots)$ is deemed positive if and only if
the set~$\{i \in \NN \,|\, x_i > 0\}$ is part of the ultrafilter.}
Non-principal ultrafilters cannot be described in explicit terms, and
they are also not at all canonical structures: \textsc{zfc} proves that there
are~$2^{2^{\aleph_0}}$ many, see~\cite{pospisil:ultrafilters}.

A practical consequence of this nonconstructivity is that it can be hard to
unwind proofs which employ hyperreal numbers to direct proofs, and even where
possible there is no general procedure for doing so. For instance, Fabert has
not obtained a direct proof of his result, and not for the lack of trying
(personal communication).

\subsection{Topos-theoretic alternatives to the hyperreal numbers} Topos theory
provides several constructive alternatives for realizing infinitesimals.
One such is ``cheap nonstandard analysis'' by~\cite{tao:cheap-nsa}.
It is to Robinson's nonstandard analysis what potential infinity is
to actual infinity: Instead of appealing to the axiom of choice to
obtain a completed ultrafilter, cheap nonstandard analysis constructs larger
and larger approximations to an ideal ultrafilter on the go.

The following section presents a (variant of a) topos used in \emph{synthetic differential
geometry} as discussed by~\cite{kock:sdg,kock:new-methods}. This subject is a further topos-theoretic
approach to infinitesimals which is suited to illustrate the philosophy of toposes
as lenses. A major motivation for the development of synthetic differential
geometry was to devise a rigorous context in which the
writings of Sophus Lie, who freely employed infinitesimals in his seminal
works, can be effortlessly interpreted, staying close to the
original and requiring no coding.

\subsection{The Zariski topos}
\label{sect:the-zariski-topos}
The starting point is the observation that
while the field~$\RR$ of ordinary real numbers does not contain infinitesimals
(except for zero), the ring~$\RR[\varepsilon]/(\varepsilon^2)$ of \emph{dual
numbers} does. This ring has the cartesian product~$\RR \times \RR$ as its
underlying set and the ring operations are defined such
that~$\varepsilon^2 = 0$, where~$\varepsilon \defeq \langle 0,1 \rangle$:
\[
  \langle a,b \rangle + \langle a',b' \rangle \defeq
  \langle a+a', b+b' \rangle \qquad
  \langle a,b \rangle \cdot \langle a',b' \rangle \defeq
  \langle aa', ab'+a'b \rangle
\]
We write~$\langle a,b \rangle$ more clearly
as~$a+b\varepsilon$.

The flavor of infinitesimal numbers supported
by~$\RR[\varepsilon]/(\varepsilon^2)$ are the \emph{nilsquare numbers}, numbers
which square to zero. The numbers~$b \varepsilon$ with~$b \in \RR$
are nilsquare in~$\RR[\varepsilon]/(\varepsilon^2)$, and
they are sufficient to rigorously reproduce derivative computations of
polynomials such as~\eqref{eq:derivative}.

However, the dual numbers are severely lacking in other aspects. Firstly, they
do not contain any nilcube numbers which are not already nilsquare. These are
required in order to extend calculations like~\eqref{eq:derivative} to second
derivatives, as in
\[ (x+\varepsilon)^3 - x^3 = 3x^2\varepsilon + \frac{1}{2!} 6x
\varepsilon^2. \]
Secondly, the dual numbers contain, up to scaling, only a single infinitesimal
number. Further independent infinitesimals are required in order to deal with
functions of several variables, as in
\[ f(x+\varepsilon,y+\varepsilon') - f(x,y) = D_xf(x,y)\varepsilon +
D_yf(x,y)\varepsilon'. \]
Thirdly, and perhaps most importantly, the ring of dual numbers fails to be a
field. The only invertible dual numbers are the numbers of the form~$a +
b\varepsilon$ with~$a$ invertible in the reals; it is not true that any nonzero
dual number is invertible.

The first deficiency could be fixed by passing
from~$\RR[\varepsilon]/(\varepsilon^2)$ to~$\RR[\varepsilon]/(\varepsilon^3)$
(a ring whose elements are triples and whose ring operations are defined such
that~$\langle0,1,0\rangle^3 = 0$) and the second by passing
from~$\RR[\varepsilon]/(\varepsilon^2)$
to~$\RR[\varepsilon,\varepsilon']/(\varepsilon^2,\varepsilon'^2,\varepsilon\varepsilon')$.
In a sense, both of these proposed replacements are better \emph{stages} than
the basic ring~$\RR[\varepsilon]/(\varepsilon^2)$ or even~$\RR$ itself. However,
similar criticisms can be mounted against any of these better stages, and the
problem that all these substitutes are not fields persists.

\bigskip
\paragraph{Introducing the topos}
The \emph{Zariski topos of~$\RR$},~$\Zar(\RR)$, meets all of these challenges. It
contains a ring~$\RR^\sim$, the so-called \emph{ring of smooth numbers}, which reifies
the real numbers, the dual numbers, the two proposed better stages and indeed
any finitely presented~$\RR$-algebra into a single coherent entity. The
Kripke--Joyal translations rules of~$\Zar(\RR)$ are listed in
Table~\ref{table:zariski}. Any evaluation of an internal statement starts out
with the most basic stage of all, the ordinary reals~$\RR$; then, during the
course of evaluation, the current stage is successively refined to better
stages (further finitely presented~$\RR$-algebras).

For instance, universal quantification~``$\forall x \? \RR^\sim$'' not only refers to
all elements of the current stage, but also to any elements of any refinement
of the current stage. Similarly, negation~``$\neg\varphi$'' does not only mean
that~$\varphi$ would imply~$\bot$ in the current stage, but also that it does
so at any later stage.

For reference purposes only, we include the precise definition of the Zariski topos.
\begin{definition}The \emph{Zariski topos of~$\RR$},~$\Zar(\RR)$, is a certain full subcategory of
the category of functors from finitely presented~$\RR$-algebras to sets, namely of
the Zariski sheaves. Such a functor is a \emph{Zariski sheaf} if and only if,
for any covering~$(A[f_i^{-1}])_i$ of any finitely presented~$\RR$-algebra~$A$
(this notion is defined in Table~\ref{table:zariski}), the diagram
\[ F(A) \rightarrow \prod_i F(A[f_i^{-1}]) \rightrightarrows \prod_{j,k}
F(A[(f_jf_k)^{-1}]) \]
is an equalizer diagram. The object~$\RR^\sim$ of~$\Zar(\RR)$ is the tautologous functor~$A \mapsto A$.
\end{definition}

\begin{table}
  \centering
  \begin{framed}
    \begin{tabbing}
      $A \models (\exists x\?\RR^\sim\_ \varphi(x))$ \= \kill
      $\Zar(\RR) \models \varphi$ \> iff $\RR \models \varphi$. \\\\
      $A \models s = t$ \> iff $s = t$ as elements of~$A$. \\
      $A \models \top$ \> iff $1 = 1$ in~$A$. \\
      $A \models \bot$ \> iff $1 = 0$ in~$A$. \\
      $A \models (\varphi \wedge \psi)$ \> iff~$A \models \varphi$ and~$A \models \psi$. \\
      $A \models (\varphi \vee \psi)$ \> iff there exists a partition~$1 = f_1 + \cdots + f_n \in A$ such that, \\
      \> \qquad for each index~$i$, $A[f_i^{-1}] \models \varphi$ or~$A[f_i^{-1}] \models \psi$. \\
      $A \models (\varphi \Rightarrow \psi)$ \> iff for any finitely
      presented~$A$-algebra~$B$, \\
      \> \qquad $B \models \varphi$ implies~$B \models \psi$. \\
      $A \models (\forall x\?\RR^\sim\_ \varphi(x))$ \> iff for any finitely
      presented~$A$-algebra~$B$ and \\
      \> \qquad any element~$x_0 \in B$, $B \models \varphi(x_0)$. \\
      $A \models (\exists x\?\RR^\sim\_ \varphi(x))$ \> iff there exists a partition~$1 = f_1 + \cdots + f_n \in A$ such that, \\
      \> \qquad for each index~$i$, there is an element~$x_0 \in A[f_i^{-1}]$ \\
      \> \qquad such that $A[f_i^{-1}] \models \varphi(x_0)$.
    \end{tabbing}

    \bigskip
    \justify
    A \emph{covering} of an~$\RR$-algebra~$A$ is a finite family of~$A$-algebras of the
    form~$(A[f_i^{-1}])_{i=1,\ldots,n}$ such that~$1 = f_1 + \cdots + f_n \in A$.
  \end{framed}
  \bigskip

  \caption{\label{table:zariski}A (fragment of) the Kripke--Joyal translation
  rules of the Zariski topos~$\Zar(\RR)$.}
\end{table}

\paragraph{Properties of the smooth numbers}\label{page:field-property}
As a concrete example, the Kripke--Joyal translation of the statement that~$\RR^\sim$ is a field,
\[ \Zar(\RR) \models \forall x\?\RR^\sim\_ (\neg(x = 0) \Rightarrow \exists y\?\RR^\sim\_ xy = 1), \]
is this:
\begin{indentblock}
For any stage~$A$ and any element~$x \in A$, \\
${\qquad}$ for any later stage~$B$ of~$A$, \\
${\qquad\qquad}$ if for any later stage~$C$ of~$B$ \\
${\qquad\qquad\qquad}$ in which~$x = 0$ holds \\
${\qquad\qquad\qquad}$ also~$1 = 0$ holds, \\
${\qquad\qquad}$ then~$B$ can be covered by later stages~$C_i$ such that, \\
${\qquad\qquad\qquad}$ for each
index~$i$, there is an element~$y \in C_i$ with~$xy = 1$ in~$C_i$.
\end{indentblock}
And indeed, this statement is true. Let a stage~$A$ (a finitely
presented~$\RR$-algebra) and an element~$x \in A$
be given. Let~$B$ be any later stage of~$A$ (any finitely
presented~$A$-algebra -- such an algebra is also finitely presented as
an~$\RR$-algebra). Assume that for any later stage~$C$ of~$B$ in which~$x = 0$
holds also~$1 = 0$ holds. Trivially,~$x = 0$ holds in the particular refinement~$C \defeq B/(x)$.
Hence~$1 = 0$ holds in~$C$. By elementary algebra, this means that~$x$ is
invertible in~$B$. Hence the conclusion holds for the singleton covering
of~$B$ given by~$C_1 \defeq B$.

\begin{remark}The Zariski topos can also be set up with an arbitrary commutative
ring~$S$ in place of~$\RR$. The resulting topos~$\Zar(S)$ contains a mirror
image~$S^\sim$ of~$S$, a reification of all finitely presented~$S$-algebras into a
single entity. The computation we just carried out also applies in this
more general context and shows that~$S^\sim$ is a field. It is in this sense that the
topos~$\Zar(S)$ provides a lens through which~$S$ looks like a field.

A small variant of this lens has been used to give a new proof of
\emph{Grothendieck's generic freeness lemma}, a fundamental theorem in
algebraic geometry about the free locus of certain sheaves. The new proof uses the lens to reduce to the case
of fields, where the claim is trivial~\cite[Section~11.5]{blechschmidt:phd},
and improves in length on all previously known proofs, even if the topos
machinery is eliminated by unrolling the appropriate
definitions as in~\cite{blechschmidt:generic-freeness}.\footnote{This contribution is
not the proper place for an exposition of Grothendieck's generic freeness
lemma, but some aspects can already be appreciated on a syntactical level.
Grothendieck's generic freeness lemma states that any finitely generated sheaf
of modules on a reduced scheme is finite locally free \emph{on a dense open}.
By employing the internal language, this statement is reduced to the following
fact of intuitionistic linear algebra: Any finitely generated module over a
field is \notnot finite free.}
\end{remark}

Within~$\Zar(\RR)$, we can construct the set~$\Delta \defeq \{ \varepsilon \? \RR^\sim \,|\,
\varepsilon^2 = 0 \}$ of nilsquare numbers. Then~$\RR^\sim$ validates the following laws:
\begin{enumerate}
\item Law of cancellation: $\forall x \? \RR^\sim\_ \forall y \? \RR^\sim\_ \bigl((\forall
\varepsilon \? \Delta\_ x\varepsilon = y\varepsilon)
\Rightarrow x = y\bigr)$
\item Axiom of micro-affinity: $\forall f \? (\RR^\sim)^\Delta\_ \exists! a \? \RR^\sim\_
\forall \varepsilon \? \Delta\_ f(\varepsilon) = f(0) + a\varepsilon$
\end{enumerate}
The unique number~$a$ in the axiom of micro-affinity deserves to be
called~``$f'(0)$''; this is how we synthetically define the derivative in
synthetic differential geometry. (However, despite these properties the Zariski
topos is not yet \emph{well-adapted} to synthetic differential geometry in the
sense of Definition~\ref{defn:well-adapted} below.)

Having motivated the Zariski topos by the desire to devise a universe with
infinitesimals, the actual ontological status of the infinitesimal numbers in
the Zariski topos is more nuanced. The law of cancellation implies that,
within~$\Zar(\RR)$, it is not the case that zero is the only nilsquare number.
However, this does not mean that there actually \emph{is} a nilsquare number
in~$\RR^\sim$. In fact, any nilsquare number cannot be nonzero, as nonzero numbers are
invertible while nilsquare numbers are not. Hence any nilsquare number in~$\RR^\sim$
is \emph{not not} zero. This state of affairs is only possible in an
intuitionistic context.

\begin{remark}The ring~$\RR^\sim$ of smooth numbers does not coincide with the Cauchy reals,
the Dedekind reals or indeed any well-known construction of the reals
within~$\Zar(\RR)$. This observation explains why~$\RR^\sim$ can satisfy the law of
cancellation even though it is an intuitionistic theorem that the only
nilsquare number in any flavor of the reals is zero.
\end{remark}

\subsection{Well-adapted models} The Zariski topos of~$\RR$ allows to compute with
infinitesimals in a satisfying manner. However it is not suited as a home for
synthetic differential geometry, a first indication being that in~$\Zar(\RR)$,
any function~$\RR^\sim \to \RR^\sim$ is a polynomial function. Hence important functions such as the
exponential function do not exist in~$\Zar(\RR)$. More comprehensively, the
Zariski topos is not a well-adapted model in the sense of the following
definition.

\begin{definition}\label{defn:well-adapted}A \emph{well-adapted model} of synthetic
differential geometry is a topos~$\E$ together with a ring~$\RR^\sim$ in~$\E$ such
that:
\begin{enumerate}
\item The ring~$\RR^\sim$ is a field.
\item The ring~$\RR^\sim$ validates the axiom of micro-affinity and several related
axioms.
\item There is a fully faithful functor~$i : \mathrm{Mnf} \to \E$ embedding the
category of smooth manifolds into~$\E$.
\item The ring~$\RR^\sim$ coincides with~$i(\RR^1)$, the image of the real line
in~$\E$.
\end{enumerate}
\end{definition}

It is the culmination of a long line of research by several authors that
several well-adapted models of synthetic differential geometry exist, see~\cite{moerdijk-reyes:models}. By the
conditions imposed in Definition~\ref{defn:well-adapted}, for any such
topos~$\E$ the following transfer principle holds: If~$f,g : \RR \to \RR$ are
smooth functions, then~$f' = g$ (in the ordinary sense of the derivative) if
and only if~$i(f)' = i(g)$ in~$\E$ (in the synthetic sense of the derivative).

Hence the nilsquare infinitesimal numbers of synthetic differential geometry
may freely be employed as a convenient fiction when computing derivatives.
Because the theorem on the existence of well-adapted models has a constructive
proof, any proof making use of these infinitesimals may mechanically be
unwound to a (longer and more complex) proof which only refers to the ordinary
reals.

\subsection{On the importance of language} The verification of the field
property of~$\RR^\sim$ in Section~\ref{sect:the-zariski-topos} on
page~\pageref{page:field-property} demonstrates a basic feature of the
Kripke--Joyal translation rules: The translation~``$\E \models \varphi$'' of a
statement~$\varphi$ is usually quite complex, even if~$\varphi$ is reasonably
transparent.

The language of toposes derives its usefulness for mathematical practice from
this complexity reduction: In some cases, the easiest way to prove a result
(about objects of the standard topos) is
\begin{enumerate}
\item to observe that the claim is
equivalent to the Kripke--Joyal translation of a different (typically more
transparent) claim about objects of some problem-specific relevant topos and
then
\item to verify this different claim, reasoning internally to the topos.
\end{enumerate}

One can always mechanically eliminate the topos machinery from such a proof, by
translating all intermediate statements following the Kripke--Joyal translation
rules and unwinding the constructive soundness proof of
Theorem~\ref{thm:reasoning}. This unwinding typically turns transparent
internal proofs into complex external proofs -- proofs which one might not have
found without the problem-adapted internal language provided by a
custom-tailored topos.

\bibliographystyle{chicago}
\bibliography{paper-filmat}

\end{document}